\documentclass[11pt]{article}
\usepackage[dvips]{graphicx}
\usepackage{textcomp}
\usepackage{amsbsy}
\usepackage{latexsym}
\usepackage[mathscr]{eucal}
\usepackage{amsfonts,amsmath,amsthm}
\usepackage{amssymb}
\usepackage[usenames]{xcolor}
\usepackage{enumitem}

\usepackage[normalem]{ulem}

\usepackage{fullpage}

\newcommand\restr[2]{{
		\left.\kern-\nulldelimiterspace 
		#1 
		\vphantom{\big|} 
		\right|_{#2} 
}}
	
\begin{document}
\title
{Nonlinear Methods for Model Reduction}
\author{ 
 Andrea Bonito, Albert Cohen, Ronald DeVore, \\
 Diane Guignard, Peter Jantsch, and Guergana Petrova
\thanks{%
 This research was supported by the NSF grants DMS 1817691 (AB), DMS 1817603 (RD-GP), ONR grant N00014-17-1-2908 (RD); DG was supported by the Swiss National Science Foundation grant P2ELP2-175056 and IAMCS at TAMU, and PJ was supported by an NSF Fellowship DMS 17-04121. A portion of this research was completed while RD (Simon Fellow), DG, and PJ were supported as visitors of the Isaac Newton Institute at Cambridge University}}%
\hbadness=10000
\vbadness=10000
\newtheorem{lemma}{Lemma}[section]
\newtheorem{prop}[lemma]{Proposition}
\newtheorem{cor}[lemma]{Corollary}
\newtheorem{theorem}[lemma]{Theorem}
\newtheorem{remark}[lemma]{Remark}
\newtheorem{example}[lemma]{Example}
\newtheorem{definition}[lemma]{Definition}
\newtheorem{proper}[lemma]{Properties}
\newtheorem{assumption}[lemma]{Assumption}
%
\newenvironment{disarray}{\everymath{\displaystyle\everymath{}}\array}{\endarray}

\def\RR{\rm \hbox{I\kern-.2em\hbox{R}}}
\def\NN{\rm \hbox{I\kern-.2em\hbox{N}}}
\def\ZZ{\rm {{\rm Z}\kern-.28em{\rm Z}}}
\def\CC{\rm \hbox{C\kern -.5em {\raise .32ex \hbox{$\scriptscriptstyle
|$}}\kern
-.22em{\raise .6ex \hbox{$\scriptscriptstyle |$}}\kern .4em}}
\def\vp{\varphi}
\def\<{\langle}
\def\>{\rangle}
\def\t{\tilde}
\def\i{\infty}
\def\e{\varepsilon}
\def\sm{\setminus}
\def\nl{\newline}
\def\wt{\widetilde}
\def\wh{\widehat}
\def\cT{{\cal T}}
\def\cA{{\cal A}}
\def\cI{{\cal I}}
\def\cV{{\cal V}}
\def\cB{{\cal B}}
\def\cF{{\cal F}}
\def\cY{{\cal Y}}

\def\cD{{\cal D}}
\def\cP{{\cal P}}
\def\cJ{{\cal J}}
\def\cM{{\cal M}}
\def\cO{{\cal O}}
\def\Chi{\raise .3ex
\hbox{\large $\chi$}} \def\vp{\varphi}
\def\lsima{\hbox{\kern -.6em\raisebox{-1ex}{$~\stackrel{\textstyle<}{\sim}~$}}\kern -.4em}
\def\lsim{\hbox{\kern -.2em\raisebox{-1ex}{$~\stackrel{\textstyle<}{\sim}~$}}\kern -.2em}
\def\[{\Bigl [}
\def\]{\Bigr ]}
\def\({\Bigl (}
\def\){\Bigr )}
\def\[{\Bigl [}
\def\]{\Bigr ]}
\def\({\Bigl (}
\def\){\Bigr )}
\def\L{\pounds}
\def\pr{{\rm Prob}}
\def\ds{\displaystyle}
\def\ev#1{\vec{#1}}   
\newcommand{\lt}{\ell^{2}(\nabla)}
\def\Supp#1{{\rm supp\,}{#1}}
\def\R{\mathbb{R}}
\def\E{\mathbb{E}}
\def\nl{\newline}
\def\T{{\relax\ifmmode I\!\!\hspace{-1pt}T\else$I\!\!\hspace{-1pt}T$\fi}}
\def\N{\mathbb{N}}
\def\Z{\mathbb{Z}}
\def\N{\mathbb{N}}
\def\Zd{\Z^d}
\def\Q{\mathbb{Q}}
\def\C{\mathbb{C}}
\def\Rd{\R^d}
\def\gsim{\mathrel{\raisebox{-4pt}{$\stackrel{\textstyle>}{\sim}$}}}
\def\sime{\raisebox{0ex}{$~\stackrel{\textstyle\sim}{=}~$}}
\def\lsim{\raisebox{-1ex}{$~\stackrel{\textstyle<}{\sim}~$}}
\def\div{\mbox{ div }}
\def\M{M} \def\NN{N}         
\def\Le{{\ell^1}}      
\def\Lz{{\ell^2}}
\def\Let{{\tilde\ell^1}}   
\def\Lzt{{\tilde\ell^2}}
\def\Ltw{\ell^\tau^w(\nabla)}
\def\t#1{\tilde{#1}}
\def\la{\lambda}
\def\La{\Lambda}
\def\ga{\gamma}
\def\BV{{\rm BV}}
\def\Ga{\eta}
\def\al{\alpha}
\def\cZ{{\cal Z}}
\def\cA{{\cal A}}
\def\cU{{\cal U}}
\def\argmin{\mathop{\rm argmin}}
\def\argmax{\mathop{\rm argmax}}
\def\prob{\mathop{\rm prob}}

\def\cO{{\cal O}}
\def\cA{{\cal A}}
\def\cC{{\cal C}}
\def\cS{{\cal F}}
\def\bu{{\bf u}}
\def\bz{{\bf z}}
\def\bZ{{\bf Z}}
\def\bI{{\bf I}}
\def\cE{{\cal E}}
\def\cD{{\cal D}}
\def\cG{{\cal G}}
\def\cI{{\cal I}}
\def\cJ{{\cal J}}
\def\cM{{\cal M}}
\def\cN{{\cal N}}
\def\cT{{\cal T}}
\def\cU{{\cal U}}
\def\cV{{\cal V}}
\def\cW{{\cal W}}
\def\cL{{\cal L}}
\def\cB{{\cal B}}
\def\cG{{\cal G}}
\def\cK{{\cal K}}
\def\cX{{\cal X}}
\def\cS{{\cal S}}
\def\cP{{\cal P}}
\def\cQ{{\cal Q}}
\def\cR{{\cal R}}
\def\cU{{\cal U}}
\def\bL{{\bf L}}
\def\bl{{\bf l}}
\def\bK{{\bf K}}
\def\bC{{\bf C}}
\def\X{X\in\{L,R\}}
\def\ph{{\varphi}}
\def\D{{\Delta}}
\def\H{{\cal H}}
\def\bM{{\bf M}}
\def\bx{{\bf x}}
\def\bj{{\bf j}}
\def\bG{{\bf G}}
\def\bP{{\bf P}}
\def\bW{{\bf W}}
\def\bT{{\bf T}}
\def\bV{{\bf V}}
\def\bv{{\bf v}}
\def\bt{{\bf t}}
\def\bz{{\bf z}}
\def\bw{{\bf w}}
\def \spn{{\rm span}}
\def \meas {{\rm meas}}
\def\rhom{{\rho^m}}
\def\diff{\hbox{\tiny $\Delta$}}
\def\EE{{\rm Exp}}
\def\lll{\langle}
\def\argmin{\mathop{\rm argmin}}
\def\codim{\mathop{\rm codim}}
\def\rank{\mathop{\rm rank}}

\def\argmax{\mathop{\rm argmax}}
\def\dJ{\nabla}
\newcommand{\ba}{{\bf a}}
\newcommand{\bb}{{\bf b}}
\newcommand{\bc}{{\bf c}}
\newcommand{\bd}{{\bf d}}
\newcommand{\bs}{{\bf s}}
\newcommand{\bff}{{\bf f}}
\newcommand{\bp}{{\bf p}}
\newcommand{\bg}{{\bf g}}
\newcommand{\by}{{\bf y}}
\newcommand{\br}{{\bf r}}
\newcommand{\be}{\begin{equation}}
\newcommand{\ee}{\end{equation}}
\newcommand{\bea}{\begin{equation*} \begin{array}{lll}}
\newcommand{\eea}{\end{array} \end{equation*}}
\def \Vol{\mathop{\rm Vol}}
\def \mes{\mathop{\rm mes}}
\def \Prob{\mathop{\rm Prob}}
\def \exp{\mathop{\rm  exp}}
\def \sign{\mathop{\rm  sign}}
\def \sp{\mathop{\rm  span}}
\def \rad{\mathop{\rm  rad}}
\def \vphi{{\varphi}}
\def \csp{\overline \mathop{\rm  span}}
%
%

\newenvironment{Proof}{\noindent{\bf Proof:}\quad}{\endproof}

\renewcommand{\theequation}{\thesection.\arabic{equation}}
\renewcommand{\thefigure}{\thesection.\arabic{figure}}

\makeatletter
\@addtoreset{equation}{section}
\makeatother

\newcommand\trunc{\mathop{\rm trunc}\nolimits}
\renewcommand\d{d}
\newcommand\dd{d}
\newcommand\diag{\mathop{\rm diag}}
\newcommand\dist{\mathop{\rm dist}}
\newcommand\diam{\mathop{\rm diam}}
\newcommand\cond{\mathop{\rm cond}\nolimits}
\newcommand\Hnorm[1]{\norm{#1}_{H^s([0,1])}}
\def\int{\intop\limits}
\renewcommand\labelenumi{(\roman{enumi})}
\newcommand\lnorm[1]{\norm{#1}_{\ell^2(\Z)}}
\newcommand\Lnorm[1]{\norm{#1}_{L_2([0,1])}}
\newcommand\LR{{L_2(\R)}}
\newcommand\LRnorm[1]{\norm{#1}_\LR}
\newcommand\Matrix[2]{\hphantom{#1}_#2#1}
\newcommand\norm[1]{\left\|#1\right\|}
\newcommand\ogauss[1]{\left\lceil#1\right\rceil}
\newcommand{\QED}{\hfill
\raisebox{-2pt}{\rule{5.6pt}{8pt}\rule{4pt}{0pt}}%
 \smallskip\par}
\newcommand\Rscalar[1]{\scalar{#1}_\R}
\newcommand\scalar[1]{\left(#1\right)}
\newcommand\Scalar[1]{\scalar{#1}_{[0,1]}}
\newcommand\Span{\mathop{\rm span}}
\newcommand\supp{\mathop{\rm supp}}
\newcommand\ugauss[1]{\left\lfloor#1\right\rfloor}
\newcommand\with{\, : \,}
\newcommand\Null{{\bf 0}}
\newcommand\bA{{\bf A}}
\newcommand\bB{{\bf B}}
\newcommand\bR{{\bf R}}
\newcommand\bD{{\bf D}}
\newcommand\bE{{\bf E}}
\newcommand\bF{{\bf F}}
\newcommand\bH{{\bf H}}
\newcommand\bU{{\bf U}}
\newcommand\cH{{\cal H}}
\newcommand\sinc{{\rm sinc}}
\def\enorm#1{| \! | \! | #1 | \! | \! |}

\newcommand{\am}{a_{\min}}
\newcommand{\aM}{a_{\max}}

\newcommand{\dm}{\frac{d-1}{d}}

\let\bm\bf
\newcommand{\bbeta}{{\mbox{\boldmath$\beta$}}}
\newcommand{\bal}{{\mbox{\boldmath$\alpha$}}}
\newcommand{\bbi}{{\bm i}}

\def\nnew{\color{Red}}
\def\mnew{\color{Blue}}
\def\wnew{\color{magenta}}

\newcommand{\dI}{\Delta}
\newcommand\aconv{\mathop{\rm absconv}}

\maketitle
\date{}
\begin{abstract} The usual approach to model reduction for parametric partial differential equations (PDEs) is to construct a linear space $V_n$ which approximates well the solution manifold $\cM$ consisting of all solutions $u(y)$ with $y$ the vector of parameters. This linear reduced model $V_n$ is then used for various tasks such as building an online forward solver for the PDE, estimating the state or parameters from data observations. It is well understood in other problems of numerical computation that nonlinear methods such as adaptive approximation, $n$-term approximation, and certain tree-based methods may provide much improved numerical efficiency. This suggests the use of nonlinear methods for model reduction as well. A nonlinear method would replace the linear space $V_n$ by a nonlinear space $\Sigma_n$. This idea has already been suggested in recent papers on model reduction \cite{eftang2010, maday2013, zou2019adaptive} where the parameter domain is decomposed into a finite number of cells and a linear space of low dimension is assigned to each cell.

Up to this point, results on such a nonlinear strategy are ad hoc and there is little known in terms of precise performance guarantees. Moreover, most numerical experiments for nonlinear model reduction have only been performed when the parameter dimension is very small (one or two). In the present paper, a step is made towards providing a more cohesive theory for nonlinear model reduction of the above type. Framing these methods in the general setting of library approximation allows us to give a first comparison of their performance with those of standard linear approximation for any general compact set.
We then turn to the study of the application of these methods for solution manifolds of parametrized elliptic PDEs.
In this context, we study a very specific example of library approximation where the parameter domain is split into a finite number $N$ of rectangular cells and where different reduced affine spaces of dimension $m$ are assigned to each cell.
The performance of this nonlinear procedure is analyzed from the viewpoint of accuracy of approximation versus $m$ and $N$. A concrete strategy for the subdivision of the parameter domain is given and it is shown how this subdivision can be exploited for various numerical tasks.  
\end{abstract}

\section{Introduction}
Complex systems are frequently described by parametrized partial differential equations (PDEs) that take the general form
\be
\cP(u,y)=0,
\ee
where $y$ ranges
over some parameter domain $Y$, and $u=u(y)$ is the corresponding solution which is assumed to be uniquely defined in some Hilbert space $V$ for every $y\in Y$. We denote by $\|\cdot\| =\|\cdot\|_V$ and $\<\cdot,\cdot\>$ the norm and inner product of $V$, respectively. In what follows, we assume that the parameters are countably infinite and have been rescaled so that $Y=[-1,1]^\N$. The case of a finite dimensional parameter $y=(y_1,\dots,y_J)$ can always be recast in this setting by considering that $u(y)$ does not depend of the variable $y_j$ for $j>J$.
 
There are three main problem areas associated with parametric PDEs: 
\begin{enumerate}[label=(\roman*)]
	\item building forward solvers to efficiently compute approximations to $u(y)$ for any given $y\in Y$;
	\item estimating the state $u(y)$ from data observation when the parameter $y$ is unknown;
	\item estimating the parameter $y$ that can give rise to an observed state $u$.
\end{enumerate}
One commonly used approach to tackle these three ranges of problems in a numerically efficient way is {\it reduced modeling}. In its most usual form, it is based on introducing a linear space $V_n$ of low dimension $n$ which is tailored to provide an accurate approximation to all solutions $u(y)$ as $y$ varies in $Y$, or equivalently, to the {\it solution manifold},
\be
\cM:=\{u(y):\ y\in Y\}.
\ee 

\subsection{Linear reduced models}  
There are two common approaches to finding a reduced model $V_n$. The first one is to establish that the forward map $y\mapsto u(y)$ has a certain analyticity in $y$, and thereby admits a Taylor series representation
\be
\label{ps}
u(y)=\sum_{\nu\in \cF} t_\nu y^\nu, \quad t_\nu\in V.
\ee 
Here $\cF $ denotes the set of $\nu=(\nu_1,\nu_2,\dots)$ which have finite support and whose entries are nonnegative integers. Quantitative bounds for the size of the Taylor coefficients $t_\nu$ allow one to prove that for each $\e$, there is a finite set $\Lambda=\Lambda(\e)\subset \cF$ such that
\be
\label{eapprox}
\sup_{y\in Y} \|u(y)-\sum_{\nu\in\Lambda} t_\nu y^\nu\|_V\le \e.
\ee
The space $V_n:={\rm span}\{t_\nu : \nu\in\Lambda\}$ provides the reduced model with $n=\#(\Lambda)$. In this case, an approximation of $u(y)$ in $V_n$ is readily provided by the function 
\be
\label{rm}
\hat u(y):= \sum_{\nu\in\Lambda} t_\nu y^\nu,
\ee
that is, using the $y^\nu$ as the coefficients of $\hat u$ in the basis $t_\nu$. Quantitative bounds on the cardinality of $\Lambda(\e)$ are known under various assumptions on the coefficients of the PDE~\cite{CD}.
 
The second approach to finding a reduced model is to judiciously select certain {\it snapshots} $u(y^1),\dots, u(y^n)$ of $u$ via a greedy procedure, and use the space $V_n:=\spn \{u(y^1),\dots, u(y^n)\}$ as the reduced model. In this case, the approximation of $u(y)$ in $V_n$ requires a projection step.
 
Recent results show that there is a numerical advantage in the Taylor coefficient approach to finding a reduced basis, at least in the case of elliptic and certain parabolic PDEs, in the sense that it is sometimes possible to a priori find the set $\Lambda$ by exploiting the parametric form of the diffusion coefficients~\cite{BCM}. This avoids computationally expensive search algorithms that are a component of greedy reduced basis selections. On the other hand, greedy procedures have the advantage that they are provably near-optimal for finding a linear space to approximate $u$, in the sense that their convergence rates are similar to those of the optimal linear spaces for approximating $\cM$, see \cite{BCDDPW}.
Moreover, as we illustrate further in this paper, numerical experiments show that for a prescribed target accuracy, the greedy generated spaces that meet this accuracy are of significantly lower dimension then their polynomial counterparts.

There is a rigorous theory that quantifies the approximation performance of both of these reduced models; see \cite{CD} for a summary of known results. The theory is most fully developed in the case of elliptic PDEs of the form 
\be 
\label{elliptic}
-\div(a\nabla u) = f,
\ee
set on a physical domain $D\subset \R^d$, say with Dirichlet boundary conditions $u_{|\partial D}=0$, and where the diffusion function $a$ has an affine parametrization
\begin{equation} 
\label{affine}
a(y)=\bar{a}+\sum_{j\geq 1}y_j\psi_j,
\end{equation}
for some given functions $\bar{a}$ and $(\psi_j)_{j\geq 1}$ in $L^{\infty}(D)$. These functions are assumed to satisfy the condition
\begin{equation} \label{hyp:UEA*}
\left\|\frac{\sum_{j\geq 1}|\psi_j|}{\bar{a}}\right\|_{L^{\infty}(D)}<1,
\end{equation}
which is equivalent to the following assumption.
\vskip .1in
\noindent
{\bf Uniform Ellipticity Assumption (UEA)}:
There exist $0<\am\leq\aM<\infty$ such that
\begin{equation} \label{hyp:UEA}
 0<\am\leq a(y)\leq \aM<\infty, \quad y\in Y.
\end{equation}
\vskip .1in
\noindent
Lax-Milgram theory then ensures that whenever $f\in V'=H^{-1}(D)$, for each $y\in Y$, the corresponding solution $u(y)$ is uniquely defined in the Hilbert space $V:=H_0^1(D)$ endowed with the norm $\|\cdot\|_V:=\|\nabla\cdot\|_{L^2(D)}$. 
 
\subsection{Nonlinear reduced models} 
It is known that in many contexts, numerical methods based on nonlinear approximation outperform linear methods, in the sense of requiring a much reduced computational cost to achieve a prescribed error tolerance \cite{Dnonlinear}. This motivates us to consider replacing the linear space $V_n$ by a nonlinear space $\Sigma_n$ depending on $n$ parameters. We call such a space $\Sigma_n$ a {\it nonlinear reduced model}. This idea has already been suggested and studied in certain settings; see e.g. \cite{eftang2010, maday2013, zou2019adaptive}. However, up till now, there has not been a unified study of nonlinear model reduction. The purpose of the present paper is to provide a formal theory for such methods of nonlinear model reduction and to prove some first results that quantify the performance of these nonlinear methods. 

The nonlinear reduced models studied in this paper can be placed into the form of what is sometimes called {\it library approximation}. Given a Banach space $X$, a library $\cL$ is a finite collection { of affine spaces} $L_1:=x_1+X_1,\dots,L_N:=x_N+X_N$, where each $X_j$ is a linear space of dimension at most $m$, and each $x_j\in X$, $j=1,\dots,N$. We set each $X_j=\{0\}$ in the case $m=0$. For an element $x\in X$, the error of approximation of $x$ by the library $\cL$ is
\be
\label{xerror}
E(x,\cL):=\inf_{L\in\cL} \dist(x,L)_X.
\ee
In other words, given $x$, we choose the best of the affine spaces $L_j=x_j+X_j$, $j=1,\ldots,N$, to approximate $x$. Given a library $\cL$ and a compact set $K\subset X$, we define
\be
\label{liberror}
E_\cL(K):=\sup_{x\in K} E(x,\cL).
\ee
Here, in the context of reduced models for parametric PDEs, the idea is to keep $m$ small when compared to the dimension $n$ used in linear models $V_n$, while retaining the same accuracy of the reduced model.
 
For parametric PDEs, we take $X=V$ and $K=\cM:=\{u(y):\ y\in Y\}$ to be the solution manifold of the PDE. A library $\cL$ would then consist of affine spaces
\be
\label{affinespace}
L_j:=u_j+V_j,
\ee
where each $u_j\in V$ and each $V_j\subset V$ has dimension at most $m$. Then, the best approximation to $u(y)$ from $L_j$ is 
\be
\label{ba}
u_j+P_{V_j}(u(y)-u_j),
\ee
where $P_{V_j}$ is the $V$-orthogonal projection onto $V_j$. In this context, when presented with a parameter $y$ for which we wish to compute an online approximation to $u(y)$, the choice of which space $L_j$ to use from a given library $\cL$ could be decided in several ways, among which we mention:
\begin{enumerate}
\item searching for a computable bound for $\dist(u(y),L_j)_V=\|u-u_j-P_{V_j}(u(y)-u_j)\|_V$, and choosing the value of $j$ that minimizes this surrogate quantity;
\item building an a priori partition of the parameter domain $Y$ into cells $Q_j$ 
and construct an $L_j$ for each cell. Then the choice of $L_j$ for approximating $u(y)$ is determined by the cell $Q_j$ containing $y$.
\end{enumerate}
Only the latter procedure is considered in this paper. 

Returning back to the case of a general Banach space $X$, we denote by $\L_{m,N}=\L_{m,N}(X)$ the collection of all libraries $\cL = \{ L_1,\ldots,L_N\}$ containing $N$ affine spaces of dimension at most $m$. If we fix the values of $m$ and $N$, then the best performance of a library with these fixed values is
\be
\label{mNwidth}
d_{m,N}(K):=\inf_{\cL \in\ \L_{m,N}}E_\cL(K).
\ee
We call $d_{m,N}$ the {\it library width} of $K$. This definition slightly differs from that introduced in \cite{Tem} in which the spaces $L_j$ are taken to be linear instead of affine.
 
Library widths include the two standard approximation concepts of widths and entropy. Recall that if $K$ is a compact set in a Banach space $X$ then its Kolmogorov $m$ width is
\be
\label{nwidth}
d_m(K):= d_m(K)_X:=\inf_{\dim(Y)=m}\dist(K,Y)_X,
\ee
where the infimum is taken over all linear spaces $Y$ of dimension $m$. Thus the Kolmogorov $m$ width of $K$ is the smallest error that can be obtained by approximation by linear spaces of dimension $m$. It follows that we can sandwich the library width $d_{m,1}(K)_X$ between Kolmogorov widths by 
\be
d_{m+1}(K)\leq d_m(K_0)=d_{m,1}(K) \leq d_{m}(K),
\ee
where $K_0=K-x_0$ for some suitable $x_0\in X$. At the other extreme,
\be
d_{0,2^n}(K)=\e_n(K),
\ee
where $\e_n(K)$ is the $n$-th entropy number of $K$ that is, the smallest number $\e$ such that $K$ can be covered by $2^n$ balls in $X$ of radius $\e$. 
 
One of the motivations for using library approximation in the context of parametric PDEs with a small value of $m$ is that the current construction of linear reduced models via greedy algorithms has offline cost that increases exponentially as the dimension of the reduced space increases. This is due to the fact that the greedy algorithm needs to search 
for the reduced basis elements through a large training set which should in principle resolve the solution manifold $\cM$ to the same accuracy $\e$ that is targeted for the reduced basis space $V_n$.
For example, it is known that if the Kolmogorov $n$ width $ d_n(\cM)$ decays like $O(n^{-s})$ for some $s>0$, then taking $\e=n^{-s}$, this training set should have cardinality $O(e^{C\e^{-1/s}})$, or equivalently $O(e^{cn})$, for some fixed constants $C,c>0$. The resulting offline cost becomes prohibitive as $\e$ is getting small (or $n$ is getting large).
The reader can find a detailed analysis of this cost of greedy constructions in \cite{CD} or \cite{CDD}. We should note that it was recently shown in \cite{CDD} that the offline cost of greedy constructions (under certain model assumptions on the parametric coefficients) can be reduced to polynomial growth in $\e$ by using random training set, provided we are now willing to accept results that hold with high probability. In order not to confuse various issues, we put this aside when going further in this paper.
 
Because of the offline cost, it may be impossible to build a linear model using a greedy algorithm when the user prescribed error is too small. On the other hand, by choosing $m$ small and an appropriate partitioning $(Q_j)_{j=1, \dots,N}$ for $Y$, the offline cost is moderate and a nonlinear reduced model may be constructed provided $N$ is not too large.
Keeping $m$ small may also be useful in other contexts such as saving in the online cost for the forward problem and numerical savings for state and parameter estimation. In fact the latter is one of our main motivations for nonlinear reduced models.
 
\subsection{Outline of the paper} 
We begin the next section by giving some general remarks on library approximation. We show that if $K$ is a compact set in a Banach space $X$ whose Kolmogorov $n$ widths decay like $n^{-r}$ for some $r>0$, then given any target accuracy $\e$ and writing $\e=n^{-r}$ for a suitable integer $n$, we have $d_{m,N}(K)\le \e$ provided $N\ge 2^{c(n-m)}$, with $c$ depending only on $r$. Thus, this result gives a bound on how many spaces would be needed in the library if we restrict the dimension of the component spaces $X_j$ to be at most $m$. While quantitative, this estimate is very pessimistic since, as is well known, nonlinear methods are not beneficial for certain compact sets.

The remainder of our paper is directed at using library approximation for reduced models for parametric PDEs. We take $K=\cM$ where $\cM$ is the solution manifold of a parametric elliptic PDE with affine coefficients \eqref{affine}. As already indicated, the library approximation studied in this paper can be viewed as first partitioning the parameter set into $N$ cells $Q_j$ and assigning an affine space $L_j=u_j+V_j$ with $V_j$ of dimension at most $m$ on each cell. The main issues therefore are how to choose the cells and how to design the spaces $V_j$. Given a target accuracy $\e$ and a prescribed target $m$ for the dimension of the spaces in the library, we are interested in strategies for generating a good partition of $Y$ into $N$ cells with a bound on the number $N$ of cells needed to guarantee the prescribed accuracy.

In \S \ref {sec:pp_PPDE}, we consider libraries where each of the $L_j$ is generated from a local polynomial expansion with $m+1$ terms. We give a tensor product strategy for subdividing the parameter domain into cells $Q_j$ which are hyperrectangles and finding a polynomial space of dimension $m+1$ associated with each cell. Thus, the reduced model can be viewed as a piecewise polynomial (in $y$) approximation to $u(y)$.
We give bounds on $N$ which are a significant improvement over those in \S \ref{sec:general} and show how these results can be used to give concrete bounds when specific assumptions are made on the affine representation \eqref{affine}. 

In \S \ref{sec:numerical}, we present the results of various numerical tests that confirm our theoretical results. First, we compare the performance (on the entire parameter domain $Y$) of the two primary linear reduced models, namely polynomial and greedy. These results show that the gain in using greedy algorithms is typically dramatic. Then we implement our numerical methods for partitioning in the case of piecewise polynomial nonlinear models, where our examples show that suitable error can be achieved with a reasonable number of cells provided $m$ is not too small. We then provide a discussion and numerical experiments of nonlinear models based on piecewise polynomials in the setting of data assimilation.

Finally in \S \ref{S:conclusions}, we conclude with remarks on the possible advantages and disadvantages of library-based reduced models for applications such as online solvers, data assimilation, and parameter estimation. This section also gives us an opportunity to mention several areas where further research is needed for a better understanding of nonlinear model reduction.

\section{General remarks on library approximation}
\label{sec:general} 
We begin by making some general remarks on library approximation. The central issue we address in this section is the size of the library needed to achieve a given target accuracy when we require dimension $m$ of the spaces in the library. The following theorem gives a first, very pessimistic, bound for the size of the library, which we denote by $N$.

\begin{theorem}
\label{libthm}
Let $K$ be a compact set in a Banach space $X$. If for some $x_0\in X$ the Kolmogorov widths of $K_0=K-x_0$ satisfy
\be
\label{lt1}
 d_k(K_0)_X \leq Mk^{-r},\quad k\geq 1,
\ee
for some $M>0$, then for any $ 0\leq m\leq n$, one has
\be
\label{lt2}
 d_{m,N}(K)\le (1+2^{2r})Mn^{-r},
\ee
provided $N\ge B_r^{n-m}$ with $B_r$ depending only on $r$. In other words, we can obtain the same accuracy as in \eqref{lt1} by using spaces of the lower dimension $m$, provided we take $N$ of them.
\end{theorem}
 
\noindent
{\bf Proof:} Since $K=K_0+x_0$ and since the definition $ d_{m,N}(K)$ uses libraries of affine spaces, it is sufficient to prove the theorem for $x_0=0$ and thus $K_0=K$.

Let us first note that there is a nested sequence of spaces $X_k\subset X_{k+1}$ with $\dim(X_k)=k$ and 
\be
\label{lt3}
\dist(K,X_k)_X\le 2^{2r}M k^{-r}, \quad   k\geq 1.
\ee
Indeed, from \eqref{lt1}, there are linear spaces $L_j$, $j\ge 0$, of dimension $2^j$, and 
\begin{equation*}
\dist(K,L_j)_X\le M2^{-jr}.
\end{equation*} 
The spaces $Y_j:=L_0+\cdots + L_j$ have dimension $n_j$ with $2^j\leq n_j\leq 2^{j+1}$ and satisfy
\be
\label{lt4}
\dist (K,Y_j)_X\le M2^{-jr} =2^{2r} M2^{-(j+2)r} \leq 2^{2r} Mn_{j+1}^{-r}, \quad j\geq 0.
\ee
Since the spaces $Y_j$ are nested, and $n_0\leq \ldots\leq n_j\leq \ldots$, we can find functions $\phi_1,\phi_2,\dots, $ such that for each $j$, the functions $\phi_1,\dots,\phi_{n_j}$ are a basis for $Y_j$. The spaces 
\begin{equation*}
X_0 := \{ 0 \},\quad X_k:={\rm span}\{\phi_1,\dots,\phi_k\}, \quad k\geq 1,
\end{equation*} 
provide such a nested sequence, since for $n_j\leq k\leq n_{j+1}$ we have $Y_j\subset X_k\subset Y_{j+1}$ and
\begin{equation*}
\dist (K,X_k)_X\leq \dist (K,Y_j)_X\leq 2^{2r} Mn_{j+1}^{-r}\leq 2^{2r} Mk^{-r}, \quad k\geq 1.
\end{equation*}

\noindent
{ \bf Case 1:} We fix $m$ and first consider the case when $n=m+2^j $ with $j=-1,0,1,\dots$, where for the purposes of this proof we replace $2^{-1}$ by $0$ when $j=-1$. We proceed by induction on $j$ and use the nested spaces $X_k$ defined above. We define $W:=X_m$ which is a space of dimension $m$ and for each $j\ge 0$, we further define
$$
Z_j:={\rm span}\{\phi_{m+1},\dots, \phi_{m+2^j}\}, \quad \dim (Z_j)=2^j, \quad \hbox{and thus}\,\,W+Z_j=X_{m+2^j}.
$$ 
 
We show by induction that for each $j\ge -1$, there is a set $S_j\subset Z_j$ such that:
\begin{enumerate}
\item the library $\cL_j:=\{ s+W,\,\,s\in S_j\}$ provides the approximation error
\be
\label{desired}
E_{\cL_j}(K)\le (1+ 2^{2r}) M[m+2^j]^{-r},\quad j\ge -1;
\ee
\item for each $j\ge -1$, the cardinality of $S_j$ is
\be
\label{card1}
\#(S_j)=:N_j\le (1+2^{r+1}R)^{2^{j+1}}, \quad R:=1+2^{2r+1}.
\ee
\end{enumerate}
When $j=-1$, we can take the set $S_{-1}=\{0\}$. We obtain the desired error bound because of \eqref{lt3} and we know that $N_{-1}=1$. 

Suppose now that we have established (i) and (ii) for a value of $j$. To advance the induction to $j+1$ we do the following. 
Let $\hat X:=X/W$ denote the quotient space of $X$ modulo $W$ with elements $[x]=x+W$, $x\in X$. We equip this space with its usual norm
\be
\label{qn}
\|[x]\|_{\hat X}:=\dist(x,W)_X.
\ee
We then have the finite dimensional spaces $\hat Z_j:=\{ [z] : z\in Z_j\}$, $j=0,1,\dots$. For each $z_\ell\in S_j\subset Z_j$, we define 
$$
B_\ell=B([z_\ell],R_0):=\{[z]\in \hat Z_{j+1}:\,\,\|[z]-[z_\ell]\|_{\hat X}\leq R_0\},\quad R_0:= RM[m+2^j]^{-r},
$$
the ball in $\hat Z_{j+1}$ with center $[z_\ell]$ and radius $R_0$. It is known (see \cite{P}, p.63) that for any $\e>0$, the covering number $N_\e(U)$ for the unit ball $U$ in $\hat Z_{j+1}$ satisfies
\be
\label{bs}
\nonumber
N_\e(U)\le (1+2/\e )^{2^{j+1}}. 
\ee
We next set $\e:=M[m +2^{j+1}]^{-r}$. It follows that the covering number of $B_\ell$ satisfies
\be
\label{coverB}
N_\e({B_\ell})\le (1+2RM[m+2^j]^{-r}/\e)^{2^{j+1}}\leq (1+2^{r+1}R)^{2^{j+1}}, \quad \ell=1,\ldots, N_j.
\ee
We now take $S_{j+1}\subset Z_{j+1}$ as a collection $\{s\}$ of representatives of the centers $[s]$ of the totality of all the balls of radius $\e$ needed to cover all of the balls $B_\ell$, $\ell=1,\ldots, N_j$, that is
$$
\bigcup_{\ell=1}^{N_j}B_\ell\subset \bigcup_{s\in S_{j+1}} B([s],\e).
$$
Clearly,
\be
\label{cardinality}
\#(S_{j+1})\le N_j (1+2^{r+1}R)^{2^{j+1}}\le (1+2^{r+1}R)^{2^{j+2}},
\ee
where we have used the induction hypothesis (ii) in the least inequality.
This advances the induction assumption for the bound on $\#(S_j)$. 

We now check that the library $\cL_{j+1}:= \{ s+W,\,\,s\in S_{j+1}\}$ provides the desired approximation error bound. Let $x\in K$. Then, it follows from \eqref{lt3}
that there is a $z\in Z_{j+1}$ such that 
\be
\label{need1}
\|[x]-[z]\|_{\hat X}\le 2^{2r} M[m+2^{j+1}]^{-r}.
\ee 
We also know from our induction hypothesis (i) that there is a 
$z_\ell\in S_j$,such that 
$$
\|[x]-[z_\ell]\|_{\hat X}\le (1+2^{2r})M[m+2^j]^{-r}.
$$
Hence,
\be
\nonumber
\|[z]-[z_\ell]\|_{\hat X}\leq \|[x]-[z]\|_{\hat X}+\|[x]-[z_\ell]\|_{\hat X}
\le (1+2^{2r+1})M[m+2^j]^{-r} ,
\ee
and so $[z]$ is in the ball $ B_\ell$. Therefore, there is an $s \in S_{j+1}$ such that 
$$
\|[z]-[s] \|_{\hat X}\le M[m+2^{j+1}]^{-r}.
$$
Combining this with \eqref{need1}, we obtain
\be
\label{obtain}
\|[x]-[s]\|_{\hat X}\le (1+2^{2r})M[m+2^{j+1}]^{-r}. 
\ee
This advances our induction hypothesis on the error bound. 

\noindent
{\bf Case 2:} We consider any $n$, not necessarily of the form $m+2^j$. For any $j$ such that $m+2^j\ge n$, the library $\cL_j$ will provide the error $ (1+2^{2r})Mn^{-r} $ because of \eqref{desired}.
So, we choose $j$ as the smallest integer such that $2^j\ge n-m$. For this value of $j$, we have $2^{j-1}\le n-m$ and from \eqref{card1}, we obtain the bound
\be
\label{st1}
N_j \le (1+2^{r+1}R)^{2^{j+1}}=B_r^{2^{j-1}} \le B_r^{n-m},
\ee
with $B_r:= (1+2^{r+1}R)^4$. \hfill $\Box$
\nl

\begin{remark}
We may restate Theorem \ref{libthm} as follows. If 
\be
\nonumber
d_k(K_0)\le Mk^{-r},\quad k\geq 1,
\ee
then for any $\e>0$ and $m\geq 0$, there exists a library $\cL$ of $m$ dimensional affine spaces which approximates $K$ to accuracy $\e$, and has cardinality
\be
\nonumber
N=\#(\cL)\leq \exp(\alpha\e^{-1/r}-\beta m),
\ee
with $\beta=\ln(B_r)$ and $\alpha=\ln(B_r)\left[M{(1+2^{2r})}\right]^{1/r}$. In particular, the library widths of $K$ satisfy
\be
\nonumber
d_{m,N}(K)\leq \e, \quad \mbox{whenever} \quad N\geq \exp(\alpha\e^{-1/r}-\beta m).
\ee
\end{remark}

Theorem \ref{libthm} is very pessimistic since it holds for all compact sets $K$ and general Banach spaces $X$. As we know in other settings, some compact model classes do not benefit from nonlinear approximation. Also, note that in the proof of the theorem, we use the same space $W$ of dimension $m$ for each of the affine spaces $L_j$, thereby never taking advantage of any local behavior of the set $K$.
In the following sections of this paper, we study library approximation for the purpose of creating a nonlinear model reduction for parametric elliptic PDEs. We exploit known theorems on the smoothness of the mapping $y\mapsto u(y)$ to give explicit non-uniform and anisotropic tensor product partitions of the parameter domain $Y$ into $N$ cells and create a library of affine spaces that achieves a prescribed target error and whose size obeys much better bounds than those given in this section.

\section{Piecewise polynomial approximation for parametric PDE} 
\label{sec:pp_PPDE}
Before beginning our analysis, we first remark on what we can expect as quantitative results. Nonlinear methods are most effective when the target function, in our case $u$, is not smooth; for example when it has point singularities or singularities on lower dimensional sets, or it is piecewise smooth. 
For the parameter to solution map $y\mapsto u(y)$ associated to the elliptic equation \eqref{elliptic} with affine parametrization \eqref{affine}, singularities occur when the function $a(y)$ is not strictly positive. The uniform ellipticity assumption \eqref{hyp:UEA} ensures that the singularities of $u$ are located outside the parameter domain $Y$. However, as $\am/\aM$ becomes small, they get closer to the boundary of $Y$, and the use of nonlinear methods becomes more relevant in those cases. 

We shall see that the bounds on the number of cells necessary in a partition generated by the nonlinear method remain modest when a reasonable number of terms $m$ in the polynomial approximation are used on each cell; see Table~\ref{table1}. In the final section of this paper, we discuss the advantages this fact provides for online solvers and state estimation.

\subsection{Polynomial approximation error}
If $\Lambda\subset \cF $ is a finite set of indices, we denote by $\cP_\Lambda$ the linear space of all $V$ valued polynomials 
\be
\label{plambda}
P(y)=\sum_{\nu\in\Lambda} c_\nu y^\nu,
\ee
where the coefficients $c_\nu$ are in $V$. Here and later we use standard multivariate notation, for example, $y^\nu= y_1^{\nu_1}\cdots$ when $\nu$ has finite support. We always assume that the set $\Lambda$ is a {\it downward closed (or lower) set, that is,} 
\be
\nu\in\Lambda \quad {\rm and} \quad \mu\le \nu \implies \mu\in\Lambda,
\ee
where $\mu\le \nu$ means that $\mu_j\leq \nu_j$ for all $j$. In particular, the null multi-index is contained in $\Lambda$. Once the coefficients $c_\nu$ are fixed, each $P(y)$ is in the affine space
\be
c_0+{\rm span}\{c_\nu \in V \, : \, \nu\in \Lambda^*\}, \quad \quad \Lambda^*:=\Lambda\setminus \{0\},
\label{affinepol}
\ee
which has dimension no more than $\#(\Lambda^*)=\#(\Lambda)-1$. A typical choice for the $c_\nu$ are the Taylor coefficients in the expansion \eqref{ps}. 

There are two types of assumptions on the diffusion coefficient commonly
employed when proving results on polynomial approximation to $u$. The first one is to assume a decay rate for the sequence $(\|\psi_j\|_{L_\infty(D)})_{j\ge 1}$.
The second type of assumption (and the one we employ here), described in \cite{BCM}, is to assume a local interaction bound on how the supports of the $\psi_j$ overlap. One could derive bounds similar to those given below in the first setting as well. 

We assume throughout this section that $u(y)$ is the solution to \eqref{elliptic} with diffusion coefficient $a(y)$ given by \eqref{affine} and that there is a positive sequence $(\rho_j)_{j\geq 1}$ such that
\be
\label{kp}
\kappa:=\min_{j\geq 1} \rho_j>1,
\ee
and 
\begin{equation} \label{eqn:WUEA}
\delta:=\left\|\frac{\sum_{j\geq 1}\rho_j|\psi_j|}{\bar{a}}\right\|_{L^{\infty}(D)}<1.
\end{equation}
The following theorem gives a bound for the error of approximation of $u$ by polynomials from $P_\Lambda$.
 
\begin{theorem}
\label{tspecific}
Assume that \eqref{kp} and \eqref{eqn:WUEA} hold with $(\rho^{-1}_j)_{j\geq 1}\in \ell_q(\mathbb{N})$ for $0 < q < 2$. For each $m\geq1$, there is a set $\Lambda$ with $\#(\Lambda)=m$ such that the $V$ valued polynomial $P(y):=\sum_{\nu\in\Lambda}t_\nu y^\nu$, $y\in Y$, satisfies
\be
\label{ts1}
\sup_{y\in Y} \|u(y)-P(y)\|_V\le C(\delta,\rho,q) \|(\rho_j^{-1})_{j\geq1}\|_{\ell_q} m^{-r}, \quad r=1/q-1/2,
\ee
where $C(\delta,\rho,q):=C(\rho,q)C_{\delta}$ with
\begin{equation} \label{eqn:cst}
C(\rho,q):=\beta^{\frac{1}{q}}\exp\(\frac{\beta}{q}\|(\rho_j^{-1})_{j\ge 1}\|_{\ell_q}^q\), \quad \beta:=-\ln(1-\kappa^{-q})\kappa^q, \quad C_{\delta}^2:=\frac{(2-\delta) \aM}{(2-2\delta) \am^3}\|f\|^2_{V'}.
\end{equation}
The set $\Lambda$ can be chosen to be a lower set and is derived explicitly in the proof. 
\end{theorem}
 
\noindent
{\bf Proof:} The proof follows from a general summability result established in \cite{BCM} together with concrete estimates for the constants given in \cite{BDGJP2019}. For the completeness and clarity of the present paper, we provide the details.
We first choose $\Lambda$ to be the set of indices in $\cF$ that correspond to the $m$ largest of the numbers $\rho^{-\nu}$. Ties are handled in such a way that $\Lambda$ is a lower set, see \cite{BDGJP2019}. Then, for $P(y):=\sum_{\nu\in\Lambda}t_\nu y^\nu$ we have by H\"older's inequality that for any $y\in Y$,
\begin{equation} \label{step1}
\|u(y)-P(y)\|_V\leq \sum_{\nu\notin\Lambda}\|t_{\nu}\|_V \leq \(\sum_{\nu\in\cF} \rho^{2\nu}\|t_\nu\|_V^2\)^{\frac{1}{2}} \(\sum_{\nu\notin\Lambda}\rho^{-2\nu}\)^{\frac{1}{2}}.
\end{equation}
From the proof of Theorem 2.2 in \cite{BCM}, we know also that
\be
\label{ts11}
\sum_{\nu\in\cF} \rho^{2\nu}\|t_\nu\|_V^2\le \frac{(2-\delta) \|\bar a\|_{L^\infty(D)}}{(2-2\delta) \inf_{x\in D}\bar a(x)^3}\|f\|^2_{V'} \leq C_{\delta}^2,
\ee
where $C_{\delta}$ is defined in \eqref{eqn:cst}.
Moreover, we have
\begin{equation} \label{step2}
\sum_{\nu\notin\Lambda}\rho^{-2\nu}=\sum_{\nu\notin\Lambda}\rho^{-\nu (2-q)}\rho^{-\nu q}\leq 
\(\sup_{\nu\notin\Lambda}\rho^{-\nu (2-q)}\)\sum_{\nu\notin\Lambda}\rho^{-\nu q}.
\end{equation}

We now let $(\gamma_k)_{k\geq1}$ be a non-increasing rearrangement of the sequence $(\rho^{-\nu})_{\nu\in\cF}$. We note that $\gamma_1=\rho^{-0}=1$ due to the fact that $\rho_1 > 1$ and $(\rho_j)_{j\geq1}$ is non-decreasing. Then we have
\begin{equation}
\sup_{\nu\notin\Lambda}\rho^{-\nu q} = \gamma^q_{m+1} \leq m^{-1} \sum_{k=2}^{m+1} \gamma_k^q \leq m^{-1} \sum_{k\geq 2} \gamma_k^q 
= m^{-1} \sum_{\nu \neq 0} \rho^{-q\nu},
\end{equation}
and hence
\begin{equation}\label{step3}
\sup_{\nu\notin\Lambda}\rho^{-\nu (2-q)} \leq \( m^{-1} \sum_{\nu \neq 0} \rho^{-q\nu} \)^{\frac{2-q}{q}}.
\end{equation}
Using \eqref{ts11} and \eqref{step3} with \eqref{step2} in \eqref{step1}, we get
\begin{eqnarray}
\|u(y)-P(y)\|_V & \leq & C_\delta \( m^{-1} \sum_{\nu \neq 0} \rho^{-q\nu} \)^{\frac{2-q}{2q}}\(\sum_{\nu\notin\Lambda}\rho^{-\nu q}\)^{\frac{1}{2}} \notag \\
& \leq & C_\delta m^{-\frac{1}{q}+\frac{1}{2}}\(\sum_{\nu\neq 0}\rho^{-\nu q}\)^{\frac{1}{q}}. \label{step4}
\end{eqnarray}
The final step of the proof is giving an upper bound of the term $\sum_{\nu\neq 0} \rho^{-q\nu}$. For this, let $\alpha:=\kappa^{-q} <1$, so that $\rho_j^{-q} \le \alpha$ for all $j\ge 1$. Now define $\beta\ge 1$ so that $1-\alpha=e^{-\beta \alpha}$, i.e., $\beta$ is the same as defined in \eqref{eqn:cst}. 
Then, $\beta$ depends only on $\kappa$, and $q$, and by the convexity of $e^{-\beta x}$, we have $1-x\ge e^{-\beta x}$ for $0\le x\le \alpha$. It follows that $(1-\rho_j^{-q})^{-1}\le e^{\beta \rho_j^{-q}}$, and therefore
\be
\label{ts14}
\sum_{\nu\neq 0} \rho^{-q\nu} = \prod_{j=1}^\infty(1-\rho_j^{-q})^{-1} - 1 \le e^{\beta b}-1\le \beta be^{\beta b},\quad b:=\|(\rho_j^{-1})_{j\geq 1}\|_{\ell_q}^q.
\ee
Taking the $q$th root in \eqref{ts14} and inserting into \eqref{step4} gives~\eqref{ts1}.
\hfill $\Box$

\begin{remark}
An important point about the above theorem is that the lower set $\Lambda$ guaranteed in the theorem can be described a priori by choosing the indices corresponding to the $n$ largest of the numbers $\rho^{-\nu}$ with ties handled properly; see also \cite{CM} and \cite{BDGJP2019}.
\end{remark}

We next want to derive a local version of the last theorem, namely we want to derive an estimate for how well the Taylor series centered at a general point $\bar y\in Y$ approximates $u$ near $\bar y$.  
Suppose that $Q_\lambda(\bar y) \subset Y$ is a hyperrectangle centered at some $\bar y\in Y$ with sidelength $2\lambda_j$ in direction $j$, i.e.,
\be
\label{def:Q}
Q_\lambda (\bar y):= \{y\in\R^\N: \quad |y_j-\bar y_j|\leq \lambda_j, \quad j\geq 1\}.
\ee
We refer to the sequence $\lambda:=(\lambda_j)_{j\ge 1}$ as the sidelength vector for this set.
 
A first local error estimate for the Taylor series at $\bar y$ is given in the following corollary. In preparation for the proof of that corollary, let us note that for $y\in Q_\lambda(\bar y)$, we have
\be
\label{onQ}
a(y)= a(\bar y)+\sum_{j=1}^\infty \frac{(y_j-\bar y_j)}{\lambda_j} (\lambda_j \psi_j)= a(\bar y)+\sum_{j=1}^\infty \tilde y_j\tilde \psi_j=:\tilde a(\tilde y),
\ee
where $\tilde y_j := \frac{y_j-\bar y_j}{\lambda_j}\in [-1,1]$ and $\tilde \psi_j:=\lambda_j\psi_j$. Therefore, 

\begin{equation*}
u(y)=\tilde u(\tilde y), \quad y \in Q_\lambda(\bar y),
\end{equation*}
with $\tilde u(\tilde y)$ the solution to
\be 
\label{prob_strong_tilde}
-\textrm{div}\left(\tilde a(\tilde y)\nabla \tilde{u}(\tilde{y})\right) = f,\quad \tilde y\in Y,
\ee
in $D$ with Dirichlet homogeneous boundary conditions.
  
We can now apply Theorem \ref{tspecific} to this new problem \eqref{prob_strong_tilde} as long as the assumptions of that theorem hold for this new problem.
\begin{cor}
\label{cor1}
Suppose the assumptions of Theorem~\ref{tspecific} hold for $\kappa$ and $\delta$ as in \eqref{kp} and \eqref{eqn:WUEA}. Consider any hyperrectangle $Q:=Q_\lambda(\bar y) \subset Y$ as in \eqref{def:Q} with center $\bar y\in Y$ and sidelength vector $\lambda$. If there is a sequence $(\tilde \rho_j)_{j\ge 1}$ (depending on $Q$) for which
{\rm \begin{enumerate}
 \label{newsequence}
	\item $\tilde \rho_j \geq \kappa$ for $j\geq 1$; 
	\item $\|(\tilde \rho^{-1}_j)_{j\geq 1}\|_{\ell_q}\leq \|( \rho^{-1}_j)_{j\geq 1}\|_{\ell_q}$;
	\item $\left \| \frac{\sum_{j\geq 1}\tilde \rho_j|\tilde \psi_j|}{a(\bar y)}\right\|_{L^\infty(D)}\le \delta, \ 
$
\end{enumerate}}
\noindent then for each $m\geq1$, there is a polynomial $P$ (depending on $Q$) with $m$ terms (whose indices are given by a lower set) such that
\be
\label{cor11}
\sup_{y\in Q}\|u(y)-P(y)\|_V\le C( \delta, \rho, q)\|(\tilde \rho^{-1}_j)_{j\ge 1}\|_{\ell_q}m^{-r}, \quad r=1/q-1/2,
\ee
where $C( \delta, \rho, q)$ is the constant from Theorem \ref{tspecific}.
\end{cor}
  
\noindent 
{\bf Proof:} This follows from Theorem \ref{tspecific} applied to the new problem \eqref{prob_strong_tilde}. We obtain the same constant because of the assumptions (i)-(iii) placed on the sequence $(\tilde \rho_j)_{ j\ge 1}$.
  
\subsection{An upper bound on the library size}
\label{SS:apriori}
We now turn to the central issue of given $m$, and a desired accuracy $\e$, how can we partition the parameter domain $Y$ into a finite number of cells such that $u$ can be approximated to this accuracy by a piecewise polynomial on this partition, where each polynomial has $m+1$ terms?
Deriving such a partition and bounding its size requires some preparatory work. Let $C:=C(\delta, \rho, q)$ be the constant of Theorem~\ref{tspecific}. We assume without loss of generality that 

\begin{equation} 
\label{eqn:cstr_eps_m_PDE}
C\|(\rho_j^{-1})_{j\ge 1}\|_{\ell_q}(m+1)^{-r}>\varepsilon,
\end{equation}
since otherwise the parameter domain $Y$ does not need to be partitioned. Namely, from Theorem~\ref{tspecific}, there is a polynomial with $m+1$ terms which approximates $u$ on $Y$ to accuracy $\e$. Since $(\rho_j^{-1})_{j\geq 1} \in \ell_q(\mathbb{N})$, we define $J:=J(\e,m) \geq 1$ to be the smallest integer such that
\begin{equation} \label{def:choice_J_PDE}
\sum_{j\ge J+1} \rho_j^{-q}\leq \frac{1}{2}C^{-q}(m+1)^{qr}\e^q.
\end{equation}
We will see that the directions $J+1,J+2,J+3,...,$ contribute at most $\e/2$ to the total error and we will not need to subdivide in these directions. For the first $J$ directions, the strategy we use distributes the remaining error equally. To that purpose, we define the quantity
\begin{equation} \label{def:sigma}
\sigma^q := \frac{1}{2J}C^{-q}(m+1)^{qr}\e^q.
\end{equation}
With this notation, we can rewrite~\eqref{eqn:cstr_eps_m_PDE} and \eqref{def:choice_J_PDE}, respectively, as
\be\label{together}
\|(\rho_j^{-1})_{j\ge 1}\|_{\ell_q}^q>2J\sigma^q, \quad\mbox{and}\quad \sum_{j\ge J+1} \rho_j^{-q}\leq J\sigma^q.
\ee
	
We begin with the following lemma.
\begin{lemma} 
\label{celllemma}
Suppose $Q\subset Y$ is a hyperrectangle with center $z=(z_1,\dots,z_J,0,0,\dots)$ and sidelength vector $\lambda=(\lambda_1,\dots,\lambda_J, 1,1,\dots)$. If
\be
\label{lambda}
\lambda_j\le \sigma(\rho_j-|z_j|) \quad j=1,\dots,J,
\ee
then here exists a $V$ valued polynomial $P_{Q}$ with $m+1$ terms such that
\be
\label{tcount1_PDE_local}
\|u(y)-P_{Q}(y)\|_V\le \e, \quad y\in Q.
\ee 
\end{lemma}	

\noindent
{\bf Proof:} We define
\begin{eqnarray}
\tilde \rho_j:=
\begin{cases}
\sigma^{-1}, \quad \hbox{if}\,\, 1\leq j\leq {J},\\
\rho_j, \quad \hbox{otherwise},
\end{cases}
\label{rhotilde}
\end{eqnarray}
and verify that $(\tilde \rho_j)_{ j\ge 1}$ satisfies the assumptions (i)-(iii) of Corollary~\ref{cor1} for $Q$.

We start with (i). It follows from the definition \eqref{kp} of $\kappa$ and from~\eqref{together} that
$$
\sigma^q < \frac{1}{2J}\|(\rho_j^{-1})_{j\ge 1}\|_{\ell_q}^q = \frac{1}{2J}\left(\sum_{j=1}^J\rho_j^{-q}+\sum_{j\geq J+1} \rho_j^{-q}\right) \leq \frac{1}{2}\kappa^{-q}+\frac{1}{2}\sigma^q,
$$
and so $\sigma^{-1} > \kappa$. Since we already know $\rho_j\ge \kappa $ for all $j$, this verifies condition (i).

We now focus on (ii). We set $\eta:=C^{-1}\varepsilon(m+1)^r$ and use the choice of $J$ in \eqref{def:choice_J_PDE} to write
\begin{equation}
\label{pp}
\|(\tilde \rho ^{-1}_j)_{j \ge 1} \|^q_{\ell_{q}}= J\sigma^{q}+\sum_{j\geq J+1}\rho_j^{-q}
\leq J\sigma^{q}+\frac{1}{2}\eta^q =\eta^q.
\end{equation}
Moreover, if we combine \eqref{pp} with \eqref{together}, we obtain 
$$
\|(\tilde \rho ^{-1}_j)_{j \ge 1} \|^q_{\ell_{q}}\leq \eta^q<\|(\rho ^{-1}_j)_{j \ge 1} \|^q_{\ell_{q}},
$$
and so (ii) holds.

Finally, to prove (iii), recall that $\tilde \psi_j=\lambda_j\psi_j$ and therefore from the inequalities \eqref{rhotilde} and \eqref{lambda} we have
$$
\tilde \rho_j|\tilde\psi_j|= \tilde \rho_j\lambda_j|\psi_j|\le (\rho_j-|z_j|) |\psi_j |.
$$
This gives
\begin{equation}\label{eqn:WUEA_tilde}
\left\|\frac{\sum_{j\geq 1}\tilde{\rho}_j|\tilde{\psi}_j|}{a(z)}\right\|_{L^{\infty}(D)} \leq \left\|\frac{\sum_{j\geq 1}\rho_j|\psi_j| - \sum_{j\geq 1}|z_j||\psi_j|}{\bar{a} - \sum_{j\geq 1} |z_j||\psi_j|}\right\|_{L^{\infty}(D)}\leq \delta. 
\end{equation}
In view of the definition of $\delta$, see \eqref{eqn:WUEA}, the last inequality  follows from 
\begin{equation*}
0 \leq \sum_{j\ge 1}| z_j||\psi_j(x)| < \sum_{j\ge 1}\rho_j|\psi_j(x)| \leq \bar a(x), \quad x\in D,
\end{equation*}
and the inequality $\left|\frac{\alpha-\beta}{\gamma-\beta}\right|\leq \left|\frac{\alpha}{\gamma}\right|$ which is valid for any $0\leq \beta<\alpha\le \gamma$. Thus, (iii) has been established.
	
We can now use Corollary \ref{cor1} to guaranteed the existence of the polynomial $P_Q$ to complete the proof. \hfill $\Box$

\vskip .1in
We are now in position to state the main theorem of this section.
    
\begin{theorem} \label{tcount_general_PDE}
Let $0<q<2$ and $(\rho^{-1}_j)_{j\geq 1}\in \ell_q(\mathbb{N})$ be a nondecreasing sequence which satisfies \eqref{kp} and \eqref{eqn:WUEA}. Let $\e>0$, $m\ge 0$ and assume that \eqref{eqn:cstr_eps_m_PDE} holds.
Then, there exists a tensor product partition of $Y$ into a collection $\cR$ of $N$ hyperrectangles such that on each $Q\in\cR$ there is a $V$ valued polynomial $P_Q$ with $m+1$ terms such that
\be
\label{tcount1_PDE}
\|u(y)-P_Q(y)\|_V\le \e, \quad y\in Q.
\ee 
Furthermore, if $J:=J(\e,m)$ is as in \eqref{def:choice_J_PDE}, then the partition is obtained by only subdividing in the first $J$ directions and the number of cells $N$ in this partition satisfies
\be
\label{tcount2_PDE}
N\le \prod_{j=1}^{ J} \left( \sigma^{-1} | \ln( 1 - \rho_j^{-1}) | + C({\sigma}) \right) \quad \mbox{for some } C({\sigma})\in (1,2).
\ee	
\end{theorem}  
\noindent {\bf Proof:} To define our tensor product grid, for each $j=1,\dots,J$, we define how we subdivide $[-1,1]$ into $(2k_j+1)$ intervals
$$I_j^{i},\quad -k_j\le i\le k_j$$
for the coordinate $y_j$. Recall that we do not subdivide any of the coordinate axis when $j>J$, i.e., $k_j=0$ and $I_j^0=[-1,1]$ when $j>J$. Also, our partition is symmetric and so $I_j^{-i}= -I_j^i$, $i=1,\dots ,k_j$.
  
We fix $j\in\{1,\dots,J\}$ and describe our partition of $[-1,1]$ into intervals corresponding to the $j$-th coordinate. Our first interval $I_j^0$ is centered at $z_j^0=0$ and has sidelength $\lambda_j^0:=\sigma\rho_j$ provided this number is less than one. Otherwise, when $\sigma\rho_j\geq 1$, we define $\lambda_j^0:=1$, and so $k_j=0$ and our partition consists only of the one interval $I_j^0=[-1,1]$. Note that since $(\rho_j)_{j\geq 1}$ is nondecreasing, when this happens it also happens for all larger values of $j$.
 
Our partition is symmetric with respect to the origin and so we only describe the intervals to the right of the origin. Our next interval $I_j^1$ has left endpoint the same as the right endpoint of $I_j^0$, has center $z_j^1$ and sidelength $\lambda_j^1$, where these numbers are defined by the relationship 
\begin{equation}\label{twoequations}
\lambda_j^1= \sigma (\rho_j - z_j^1) .
\end{equation}
The only exception to this definition is when the right endpoint of this interval is larger than one. Then we recenter the interval so its left endpoint is as before and its right endpoint is one. In this case, we would stop the process and $k_j$ would be one.
 
We continue in this way moving to the right. So, in general, the interval $I_j^i$ will have its left endpoint equal to the right endpoint of $I_j^{i-1}$, and will have center $z_j^i$ and sidelength $\lambda_j^i$ which satisfy
\be
\label{satisfies}
\lambda_j^i=\sigma(\rho_j-z_j^i)
\ee
except in the case that such a choice would give a right endpoint larger than one in which we rescale. It follows that the interval $I_j^i$ always satisfies
\be
\label{satisfies1}
\lambda_j^i\le \sigma (\rho_j - z_j^i),\quad i=0,1,\dots,k_j,
\ee 
with equality except for possibly the last interval $I_j^{k_j}$. We give below a bound for $k_j$ which shows this process is finite.

This partitioning gives a tensor product set $\cR$ of hyperrectangles $Q$. In view of the property \eqref{satisfies1}, each of the hyperrectangles satisfies the conditions of Lemma~\ref{celllemma} and therefore the existence of the polynomials $P_Q$, $Q\in\cR$ satisfying the approximation estimate is guaranteed. 

It remains to bound the cardinality of $\cR$. For this, we bound $k_j$, $1\le j \le J$, when $k_j\neq 0$. We obtain the bound we want by monitoring the points
\be
\nonumber
R^i=z_j^i+\lambda_j^i, \quad i= 0,1,\dots,k_j, 
\ee
Each $R^i$ is the right endpoint of $I_j^i$ as long as $0\le i<k_j$. Also we know that $R^{k_j}\ge 1$. Relation \eqref{satisfies} implies that $\lambda_j^i$ is chosen so that
\be
\nonumber
\frac {\lambda_j^i}{\rho_j- R^i+\lambda_j^i} =\sigma. 
\ee
This gives that
\be
\nonumber
(1-\sigma)\lambda_j^i= \sigma(\rho_j-R^i).
\ee
Since $R^i=R^{i-1}+2\lambda_j^i$, we have 
\be
\nonumber
(1-\sigma)(R^i-R^{i-1})= 2 \sigma(\rho_j-R^i).
\ee
We therefore obtain the recursive formula
\be
\nonumber
R^i=\frac{1-\sigma}{1+\sigma} R^{i-1}+\frac{2\sigma}{1+\sigma} \rho_j =:\alpha R^{i-1}+b,  \quad i=1,2,\dots,
\ee
where $R^0=\rho_j \sigma$, $\alpha:=\frac{1-\sigma}{1+\sigma}$, $b:=\frac{2\sigma}{1+\sigma} \rho_j $.
Therefore, we find
\begin{eqnarray}
\nonumber
R^i&=& \alpha^iR^0+(1+\alpha+\ldots+\alpha^{i-1})b=\alpha^iR^0+\frac{1-\alpha^i}{1-\alpha} b\\
\label{lastformula}
&=&\alpha^i\rho_j \sigma + (1-\alpha^i)\rho_j =\rho_j(1-\alpha^i(1-\sigma)).
\end{eqnarray}
The iteration will stop at the smallest integer $k=k_j$ such that $R^k\geq 1$. 
Since $\sigma^{-1}\geq \kappa>1$, we have $\sigma<1$ and the iteration will stop at the smallest integer $k$ such that
$$
\alpha^k\leq \frac{1-\rho^{-1}_j}{1-\sigma}.
$$
Note that $\frac{1-\rho^{-1}_j}{1-\sigma}<1$ because $\sigma\rho_j<1$ (otherwise $k_j=0$ and $I_j^0=[-1,1]$). We are looking for the smallest integer $k$ for which 
$$
k\geq \frac{\ln\left(1-\rho^{-1}_j\right)-\ln\left(1-\sigma\right)}{\ln \alpha},
$$
which gives
\be
\nonumber
k_j=\left\lceil \frac{ \ln (1-\rho_j^{-1} ){ -\ln\left(1-\sigma\right)}}{\ln \alpha} \right\rceil < \frac{ \ln (1-\rho_j^{-1} )-\ln\left(1-\sigma\right)}{\ln \alpha}+1, \quad j=1,\dots,{J}.
\ee
Therefore, we have the bound
\begin{equation*}
n_j := 2k_j + 1 \leq 2 \frac{\ln\left( 1 - \rho_j^{-1} \right) -\ln\left(1-\sigma\right)}{\ln\left(\frac{1-\sigma}{1+\sigma}\right)} +3 = 2 \frac{\ln\left( 1 - \rho_j^{-1} \right)}{\ln\left(\frac{1-\sigma}{1+\sigma}\right)}+C({\sigma}),
\end{equation*}
where
\begin{equation} \label{def:C_sigma}
C({\sigma}):= -2\frac{\ln\left(1-\sigma\right)}{\ln\left(\frac{1-\sigma}{1+\sigma}\right)} +3 = \frac{\ln\left(\frac{(1-\sigma)}{(1+\sigma)^3}\right) }{\ln\left(\frac{1-\sigma}{1+\sigma}\right)}\in (1,2).
\end{equation}

Since $\ln (1+x)\geq \frac{2x}{2+x}$ for $x\geq 0$, we obtain
\be
\nonumber
\ln\left (\frac{1+ \sigma}{1- \sigma}\right) = \ln\left(1 + \frac{2\sigma}{1- \sigma}\right) \geq 2 \sigma,
\ee
and thus $ n_j\leq \sigma^{-1}|\ln\left( 1 - \rho_j^{-1} \right)|+C(\sigma)$, which brings us to the final calculation

\begin{equation}\label{bound_N_proof}
N = \prod_{j=1}^{J} n_j \leq \prod_{j=1}^{J}\left(\sigma^{-1}|\ln\left( 1 - \rho_j^{-1} \right)|+{C({\sigma})}\right),
\end{equation}
which completes the proof. \hfill $\Box$

\begin{remark}\label{rem:betterbound}
It follows from the proof of Theorem~\ref{tcount_general_PDE} that a more precise estimate for the number of cells is
$$
N\leq\prod_{j=1}^{J_0}\left(\sigma^{-1}|\ln\left( 1 - \rho_j^{-1} \right)|+{C({\sigma})}\right),
$$ 
where $1\leq J_0 \leq J$ is the largest integer such that $\sigma \rho_{J_0} <1$. This comes from the fact that $k_j=0$ for $J_0<j\leq J$, i.e., we do not subdivide in the directions $J_0+1.\ldots,J$.
\end{remark}

\bigskip
Let us reformulate the above result in terms of library widths. As we have remarked earlier (see \eqref{affinepol}), a polynomial approximation with $m+1$ terms is naturally associated with an affine space of dimension at most $m$. We then obtain a library $\cL=\cup_{i=1}^N L_i$ of affine spaces $L_i=L_i(P_i,Q_i)$,
\be
\nonumber
L_i=c_0^i+{\rm span}\{c^i_\nu \in V \; : \; \nu\in \Lambda_i, \,\,\#(\Lambda_i)\leq m \}, \quad i=1,\dots,N,
\ee
each associated with the polynomial $P_i$ over a hyperrectangle $Q_i\subset Y$,
\be
\nonumber
P_i(y)= c_0^i+\sum_{\nu\in \Lambda_i} c_\nu^i y^\nu,\quad y\in Q_i,
\ee
and cardinality 
$$
N\leq \prod_{j=1}^{J} \left( \sigma^{-1} | \ln\left( 1 - \rho_j^{-1} \right) | +{C(\sigma)} \right), \quad C(\sigma)\in (1,2).
$$
Moreover, since $\sup_{y\in Q_i}\|u(y)-P_i(y)\|_V\leq \e$ for $i=1,\ldots,N$, we have
$$
E_{\cL}(\cM) =\max_{y\in Y}\min_{L\in \cL}{\rm dist}(u(y),L)_V \leq \e,
$$
and therefore
$$
d_{m,k}(\cM)\leq \e, \quad \hbox{whenever} \quad k\geq\prod_{j=1}^{J} \left( \sigma^{-1} | \ln\left( 1 - \rho_j^{-1} \right) | +C(\sigma) \right).
$$

\subsection {Examples} 
\label{sec:examples}

To see how how the bounds for $N$ in Theorem~\ref{tcount_general_PDE} grow with decreasing $\e$, we consider the following 
 standard example:
 \begin{equation} \label{def:rho_decay}
\rho_j=M j^{s}, \quad j\geq 1,
\end{equation}
where $s>1/2$ is fixed.
 From our overriding assumption that $\kappa=\rho_1>1$, it follows that $M>1$.
  We note at the outset that a similar analysis can be done for other growth assumptions on the sequence $(\rho_j)_{j\ge 1}$, e.g., $\rho_j=1+Mj^s$ with $M>0$.

Before beginning our analysis, we wish to orient the reader to what type of results we can expect by reflecting on the corresponding results for polynomial approximation. In that case, we know that for each $r<s-1/2$ we can find $V$ valued polynomials
$P_n$ with $n$ terms that satisfy
\be
\label{polerror}
 \max_{y\in Y}\|u(y)-P_n(y)\|_{V}\le C_rn^{-r},\quad n=1,2,\dots .
\ee
This follows from Theorem \ref{tspecific}  by choosing a value of $q\in (1/s,2)$ with $r=1/q-1/2$.  However, we cannot take $r=s-1/2$ since the constants $C_r$ tend to infinity as $q\to1/s$. If we are given a target accuracy $\e$ then we would find the minimal number of terms $n$ to reach this accuracy by optimizing over the choice of $q$. This type of analysis is subtle and done in \cite{BDGJP2019}. 

We shall obtain similar results for piecewise
polynomial approximation where now the main new ingredient is to bound the number of cells that are needed.
We fix the desired target accuracy $\e>0$ and the value $m$ and use the a priori bound of Theorem \ref{tcount_general_PDE} to see how many 
hyperrectangles $N$ are needed to guarantee the accuracy $\e$ using piecewise polynomials with $m+1$ terms to approximate $u$ on each rectangle.  We can apply  Theorem \ref{tcount_general_PDE}  for any $q$ that satisfies $1/s<q<2$. We consider any such $q$, fix it for the moment, and investigate the size of $N$ needed to achieve the accuracy $\e$. Throughout the derivation, we let $C$ denote a constant that depends only on $q$ and may change from line to line. Note that $C_0:=C(\delta,\rho,q)$ depends only on $q$ since $\rho$ and $\delta$ are fixed. 

Since we have 
\begin{equation}
\nonumber
\sum_{j\ge J+1}\rho_j^{-q}=M^{-q}\sum_{j\ge J+1}j^{-sq}\leq C J^{1-sq},
\end{equation}
 the condition \eqref{def:choice_J_PDE} is satisfied if
 \begin{equation} \label{def:J_specific}
J=C\left(\varepsilon (m+1)^r\right)^{\frac{q}{1-sq}} =C\lambda^{\frac{q}{1-sq}}, 
\end{equation}
where
\be
\nonumber
\lambda:=\e (m+1)^{ r}, \quad \ r=r(q):=\frac{1}{q}-\frac{1}{2}.
\ee
Defining $J$ by \eqref{def:J_specific} gives that the value of $\sigma$ in the theorem is 
\be
\label{sigma}
\sigma = 2^{-1/q} C_0^{-1}J^{-1/q}\lambda= CJ^{-s}.
\ee
 Theorem \ref{tcount_general_PDE} says that we obtain a partition into $N$ hyperrectangular cells such that there is a polynomial with $m+1$ terms on each cell which achieves the desired accuracy $\e$.  It also gives that the number $N=N(q)$ of these cells can be bounded by 
\be
\label{boundNq}
N \leq \prod_{j=1}^{ J} \left( \sigma^{-1} | \ln\left( 1 - \rho_j^{-1} \right) | +{ C(\sigma)} \right)<
{ \prod_{j=1}^J \left( \sigma^{-1} | \ln\left( 1 - \rho_j^{-1} \right) | +2 \right)}.
	\ee
Since each $\rho_j\ge M>1$, and $|\ln (1-x)|\le \frac{x}{1-x}$, for $0<x<1$, we have
   \be
   \label{bdlog} |\ln (1-\rho_j^{-1})| \le (Mj^s-1)^{-1} \le (M-1)^{-1}j^{-s}, \quad j=1,2, \dots.
   \ee
Placing this into \eqref{boundNq} gives
\be
\label{boundNq1}
N \leq \prod_{j=1}^ J \left( (M-1)^{-1} \sigma^{-1}j^{-s}+2 \right) = \prod_{j=1}^J \left( CJ^s j^{-s}+2\right)
\le C^J J^{sJ}[J!]^{-s}\le e^{(C+s)J}=e^{C\lambda^{\frac{q}{1-sq}}},
	\ee 
where the last inequality uses Stirling's formula.

We examine what this bound guarantees for different values of $m$:

\noindent
{\bf Case $m=0$}:  In this case, we are providing the solution manifold $\cM$ with an $\e$ approximation net with $N$ elements. Since $\lambda=\e$ in this case, the bound \eqref{boundNq1} says we can achieve approximation accuracy $\e$ with such a net with
 \begin{equation*}N\le\exp \left \{ C\e^{ -\frac{1}{s-1/q}}\right\}\end{equation*}
   elements for any  $q\in (1/s,2)$. The best choice of $q$ in this case is to choose
$q$ as close to $2$ as possible thereby getting $N\le e^{C\e^{-1/\alpha}}$ for any $0<\alpha<s-1/2$. Notice that this is in complete agreement with what we know about the entropy of
the solution manifold $\cM$. Indeed, from Theorem \ref{tspecific}, we know the Kolmogorov width of $\cM$ satisfies
\be
\label{Kw}
 d_n(\cM)\le C_rMn^{-r}, \quad 0<r<s-1/2,
 \ee
 where the constants $C_r$ tend to infinity as $r$ gets closer to $s-1/2$. From Carl's inequality we obtain that the $\e$ covering number of $\cM$ is bounded by $e^{C\e^{-1/r}}$ provided that $r<s-1/2$ which is exactly what the above
 bound on $N$ gives.
 
 \noindent
 {\bf Case of general $m$:}  In this case, the partitioning gives a library of $N$ affine spaces of dimension $m$ that approximate $\cM$ to accuracy $\e$. In order to compare our results on piecewise polynomial approximation with those for polynomial approximation, we suppose a value of $q\in(1/s,2)$ has been chosen which gives the accuracy $C_r n^{-r}$, $r= r(q)=1/q-1/2$ using polynomials. We obtain the same accuracy $\e:=C_r n^{-r}$ using piecewise polynomial with $m+1$ terms and  the above estimate says we can do this with 
 \be
 \label{genm}
 \nonumber
 N\le \exp \left\{C\left(\frac{n}{m+1}\right)^{\frac{r}{s-1/q}    }\right\}=  \exp \left\{C\left(\frac{n}{m+1}\right)^{\alpha    }\right\},\quad \alpha:= \frac{1/q-1/2}{s-1/q},
 \ee
cells chosen as in Theorem \ref{tcount_general_PDE}. In this estimate, notice that rather than the bound $e^{C(n-m)}$ derived in \S 2 for general libraries, we now have the bound $e^{C(n/m)^\alpha}$ which gets more favorable as $m$ gets large. Note that we can always get $\alpha=1$ by taking $q=\frac{4}{2s+1}$, which belongs to the prescribed range $(1/s,2)$, since $s>1/2$ by assumption. Moreover, $\alpha$ tends to infinity as $q\rightarrow 1/s$ and to $0$ as $q\rightarrow 2$.

\section{Numerical examples}
\label{sec:numerical}
	
In this section, we present numerical examples to illustrate the performance of the strategy described above for constructing nonlinear reduced models based on 
partitioning of the parameter domain $Y$ and using piecewise $V$ valued polynomials subordinate to the chosen partition. For our numerical tests, we consider the elliptic equations \eqref{elliptic} on the domain $D=[0,1]^2$ with right-hand side $f=1$ and an affine diffusion of the form
\be
\label{affine1}
a(x,y):=1+ \sum_{j=1}^{64} y_jc_j\chi_{D_j}(x),
\ee
where $ (D_j)_{j=1}^{64}$ is a partition of $D$ into 64 square cells of equal size. The indexing is assigned randomly and has little effect on the numerical results.  Thus, the parameter domain $Y=[-1,1]^{64}$.
	
We carry out numerical experiments for different sequences $(c_j)_{j=1,\dots,64}$ 
that depend on the parameters $a_{\rm min}$ and $s$, namely
\be
\label{oursequences}
 c_j = (1-a_{\rm min})j^{-s},\quad j=1,2,\dots, 64,
\ee
where $s\in \{2,3,4\}$ and $a_{\rm min}\in \{0.1,0.05,0.01\}$. Notice that $a_{\rm min}$ is the true minimum of $a$ on $D\times Y$. Given this sequence, we can take 
\be
\label{rhoj}
\rho_j:=\frac{1-a_{\rm min}/2}{1-a_{\rm min}}j^s, \quad j=1,2,\dots,64,
\ee
and this gives  $\delta=1-\frac{a_{\rm min}}{2}$ in \eqref{eqn:WUEA}. A small value for $a_{\rm min}$ corresponds to a reduction in the domain of analyticity of $u(y)$ near the face $y_1=-1$. So, each numerical experiment corresponds to an assignment of $a_{\rm min}$ and $s$.

\subsection{Linear reduced models}
\label{SS:numericallinear}
We begin this section by considering linear reduced models with the goal of understanding how large the dimension of the linear space has to be in order to guarantee a prescribed error $\e$. We are also interested to see the effect of different choices for the linear space.  In all of our numerical experiments we take the target error to be
\be
\nonumber
\e:=10^{-4}.
\ee

We consider two choices of linear reduced models:
\begin{itemize}
\item Taylor polynomial space; 
\item reduced basis space based on greedily selected snapshots. 
\end{itemize}

We compare the approximations obtained using a Taylor polynomial with $n$ terms and a reduced basis space of dimension $n$. In particular, we want to see how large $n$ has to be to achieve the target accuracy $\e$ for these two choices.

In the case of a Taylor  polynomial space, the approximant $\bar u_n$ is given by
\be
 \label{apr}
 \bar u_n(y):=\bar t_{0}+\sum_{\nu\in \Lambda_n^*}\bar t_\nu y^\nu \in \bar t_{0} +V_{n-1}(T), \quad V_{n-1}(T):= {\rm span}\{\bar t_\nu: \, \nu\in\Lambda_n^*\},
\ee
where $\bar t_{\nu}$ is the approximation of $t_{\nu}$ obtained using a finite element solver of high accuracy (much higher accuracy than the target accuracy $\e$).
We consider two methods to generate the  lower set $\Lambda^*_n$ of cardinality $n-1$ which gives the indices $\nu$ in \eqref{apr}.

The first method, which we refer to as the {\it a priori method},  orders the $\rho^{-\nu}$, $\nu\in\cF$, in decreasing 
order according to their size. So $\nu^0:=0$ is the index giving the largest of these numbers, and $\nu^1,\nu^2,\dots$ denote the indices corresponding to the next largest of the $\rho^{-\nu}$. Ties are handled in such a way that $\Lambda_n:=\{ \nu^{0},\nu^{1},\dots,\nu^{n-1}\}$ is a lower set, see \cite{BDGJP2019}.
 We then take $\Lambda_n^*:=\Lambda_n\setminus\{\nu^0\}$.

 In the second method, here referred to as the {\it adaptive method},  we use the so-called Algorithm LN (largest neighbor) described in \cite{CCDS}
to generate an index set $ \tilde\Lambda_n$. It begins with $\nu^0:=0$ and $ \tilde\Lambda_0:=\{\nu^0\}$. Then, for $k=0,1,\ldots,n-1$,
\begin{equation} \label{eqn:adapt}
 \tilde\Lambda_{k+1}:= \tilde\Lambda_k \cup\{\nu^k\}, \quad \mbox{where} \quad \nu^k\in \argmax_{\nu\in\cR_{{ \tilde \Lambda_k}}}\|\bar t_{\nu}\|_V.
\end{equation}
Here, $\cR_{{ \tilde \Lambda_k}}$ denotes the reduced margin of the current lower set $\tilde\Lambda_k$, namely
\begin{equation*}
\cR_{\tilde \Lambda_k}:=\{\nu\in\cF\setminus\tilde \Lambda_k: \,\, \nu-e_j\in \tilde \Lambda_k \quad \mbox{for all } j \mbox{ with } \nu_j>0 \}.
\end{equation*}
 We then take $\Lambda_n^*:=\tilde \Lambda_n\setminus\{\nu^0\}$.

We compute the error $\epsilon_n$ for each of these choices by taking a large number of random (with respect to the uniform distribution) choices\footnote {In the experiments given the number of random selections of $y$ was $10^3$ and using the Mersenne Twister pseudo random generator with seed value 515.} of parameters $y\in Y$, as follows. For each choice $y$, we take an accurate finite element approximation $\bar u(y)$ of $u(y)$ as truth.  
 Note that because $\Lambda^*_n\cup \{0\}$ is a lower set, the  Taylor coefficients  $t_\nu$, $\nu\in\Lambda^*_n\cup \{0\}$, can be found recursively, see equations (3.1) and (3.2) in \cite{CCDS}.
We calculate
 $\|\bar u(y)-\bar u_n(y)\|_V$ and the error $\epsilon_n$ is then computed by maximizing $\| \bar u(y)-\bar u_n(y)\|_V$ over the random choices of $y$.

\begin{figure}[htbp]
	\centering
	\hspace*{-0.1cm}\includegraphics [width=.35\textwidth]{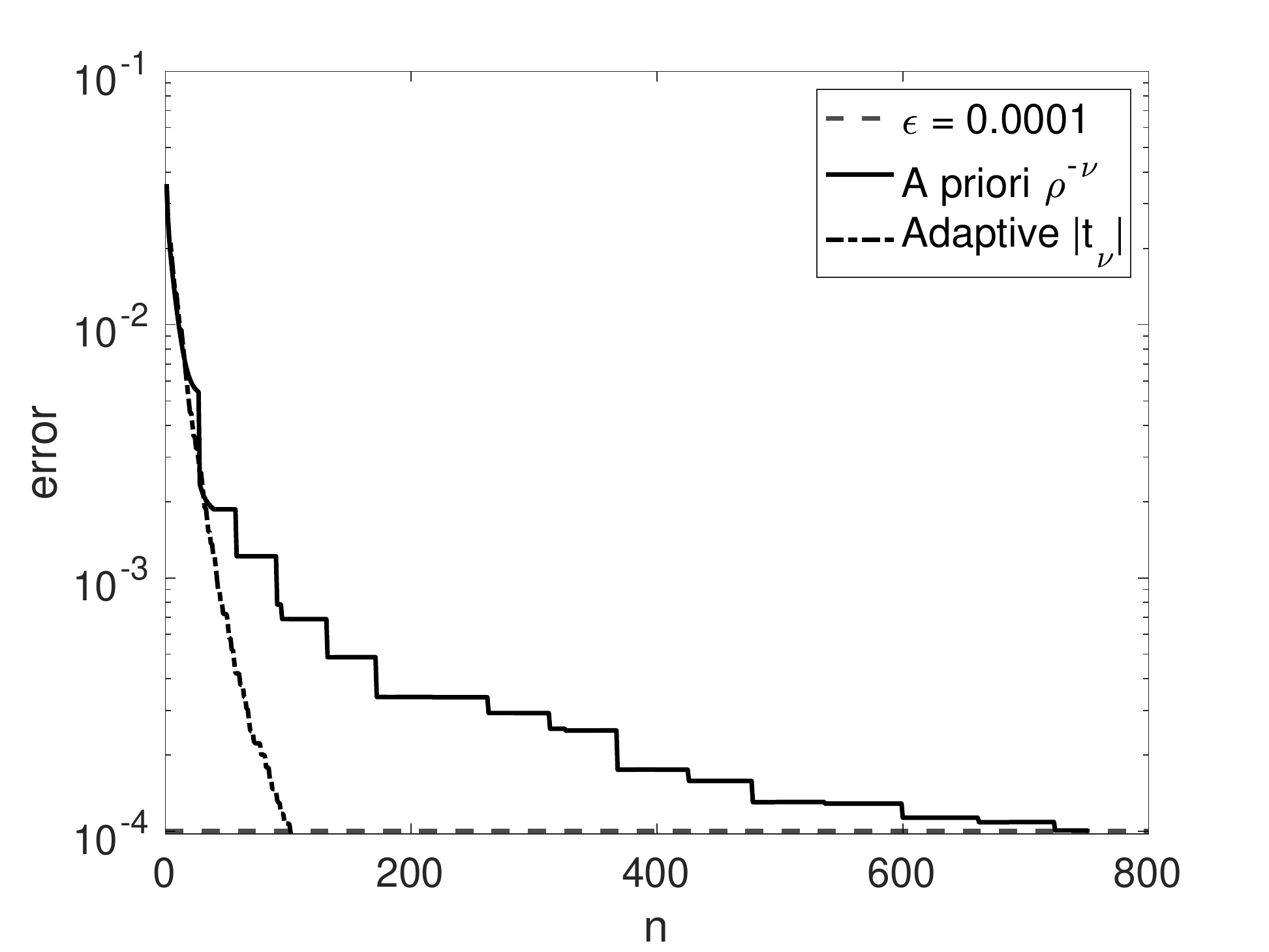}
	\hspace*{-0.5cm}\includegraphics [width=.35\textwidth]{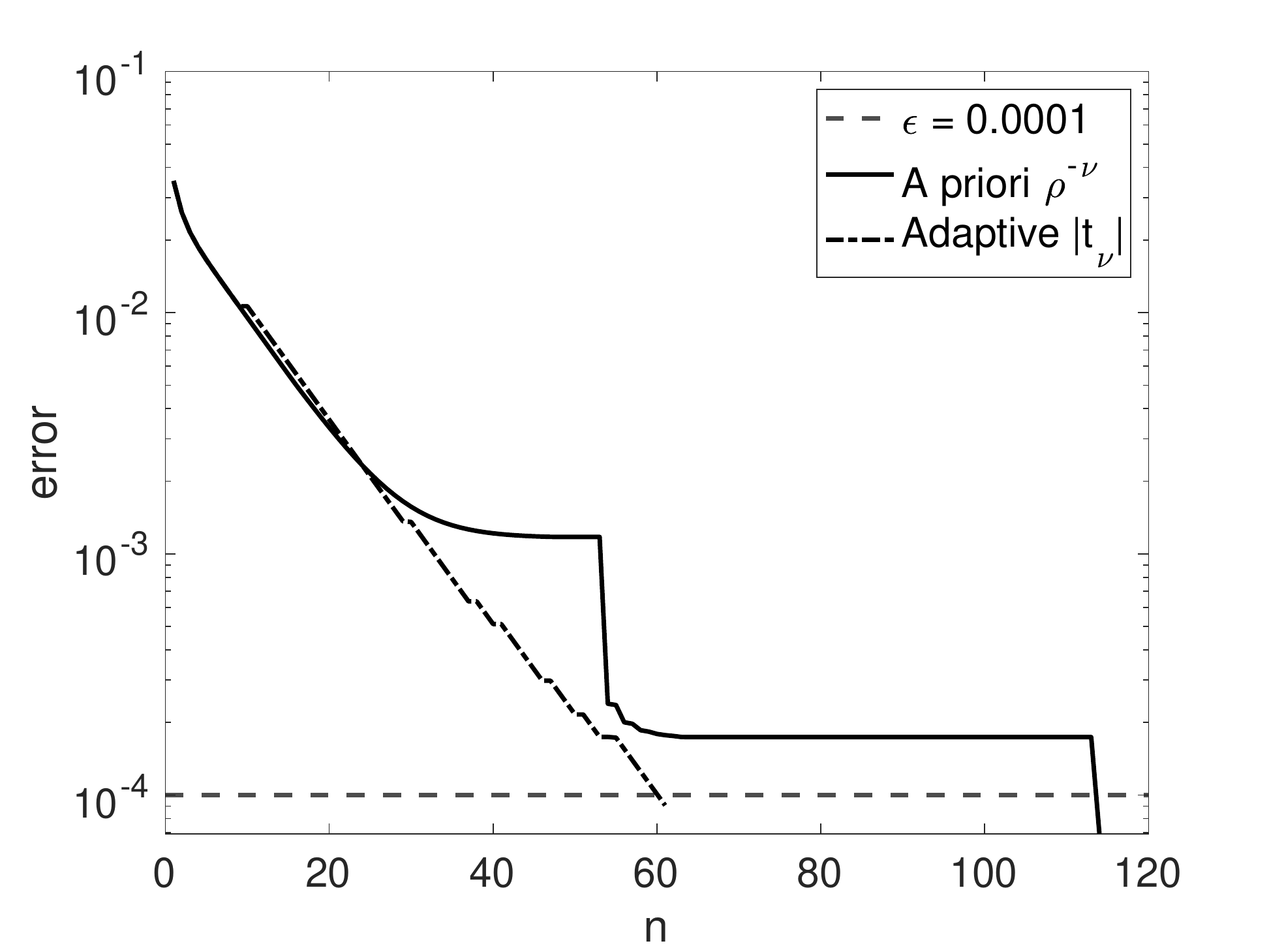}
	\hspace*{-0.5cm}\includegraphics [width=.35\textwidth]{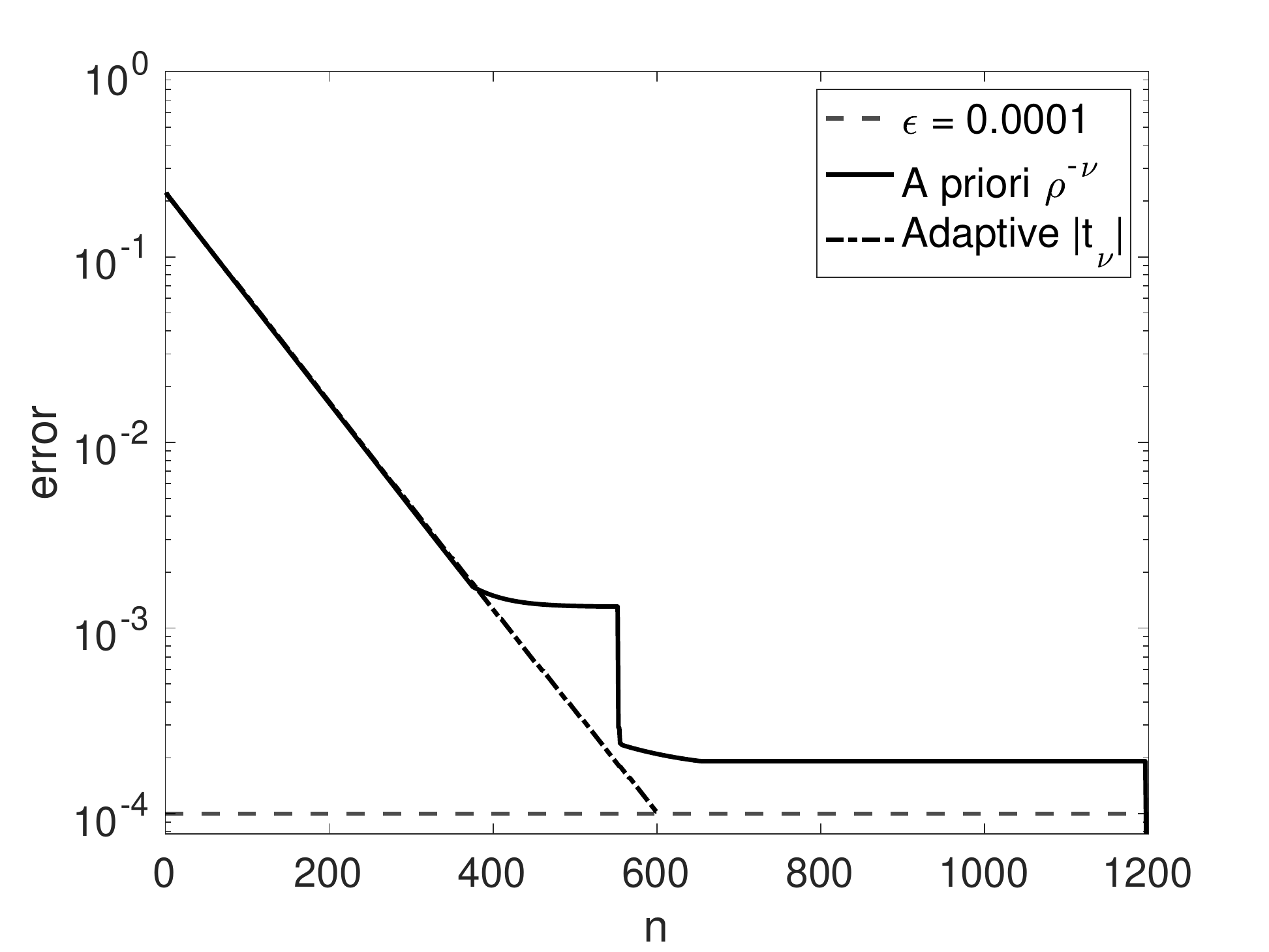} 
	\caption{Error between $\bar u$  and  the Taylor polynomial approximation $\bar u_n$ versus the number of terms $n$ for both the {\it a priori method} and the {\it adaptive method} for constructing $\Lambda^*_n$. Left: $s=2$, $a_{\rm min}=0.1$; middle: $s=4$, $a_{\rm min}=0.1$; right: $s=4$, $a_{\rm min}=0.01$.} \label{fig:apriori_adapt}
\end{figure}

Figure \ref{fig:apriori_adapt} shows a comparison of the errors obtained using the adaptive and the a priori 
methods to compute the set $\Lambda^*_n$ as $n$ grows for different values of $s$ and $a_{\rm min}$.
We see that the adaptive method to generate $\Lambda^*_n$ outperforms the a a priori method, in that the corresponding approximation error is smaller for the adaptive method. This is caused by the fact that $\|\bar t_\nu\|_V$ could be much smaller than  $\rho^{-\nu}$.  On the other hand,
the computational cost to find $\Lambda^*_n$ is greater for the adaptive method. In going further in this section, we always compute the set $\Lambda^*_n$ for Taylor polynomial indices by using the adaptive method. 

\begin{figure}[htbp]
	\centering
	\hspace*{-0.1cm}\includegraphics [width=.35\textwidth]{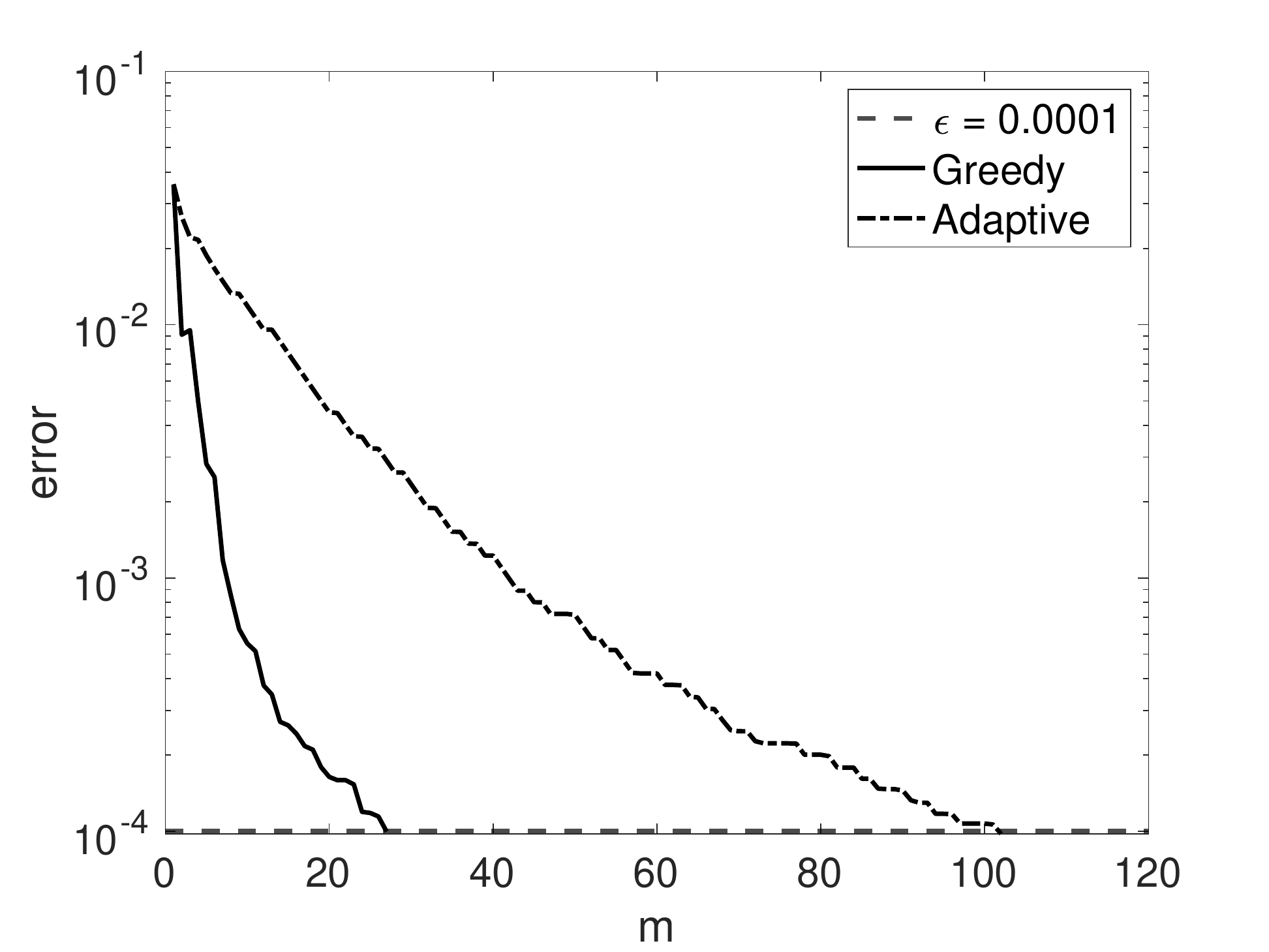}
	\hspace*{-0.5cm}\includegraphics [width=.35\textwidth]{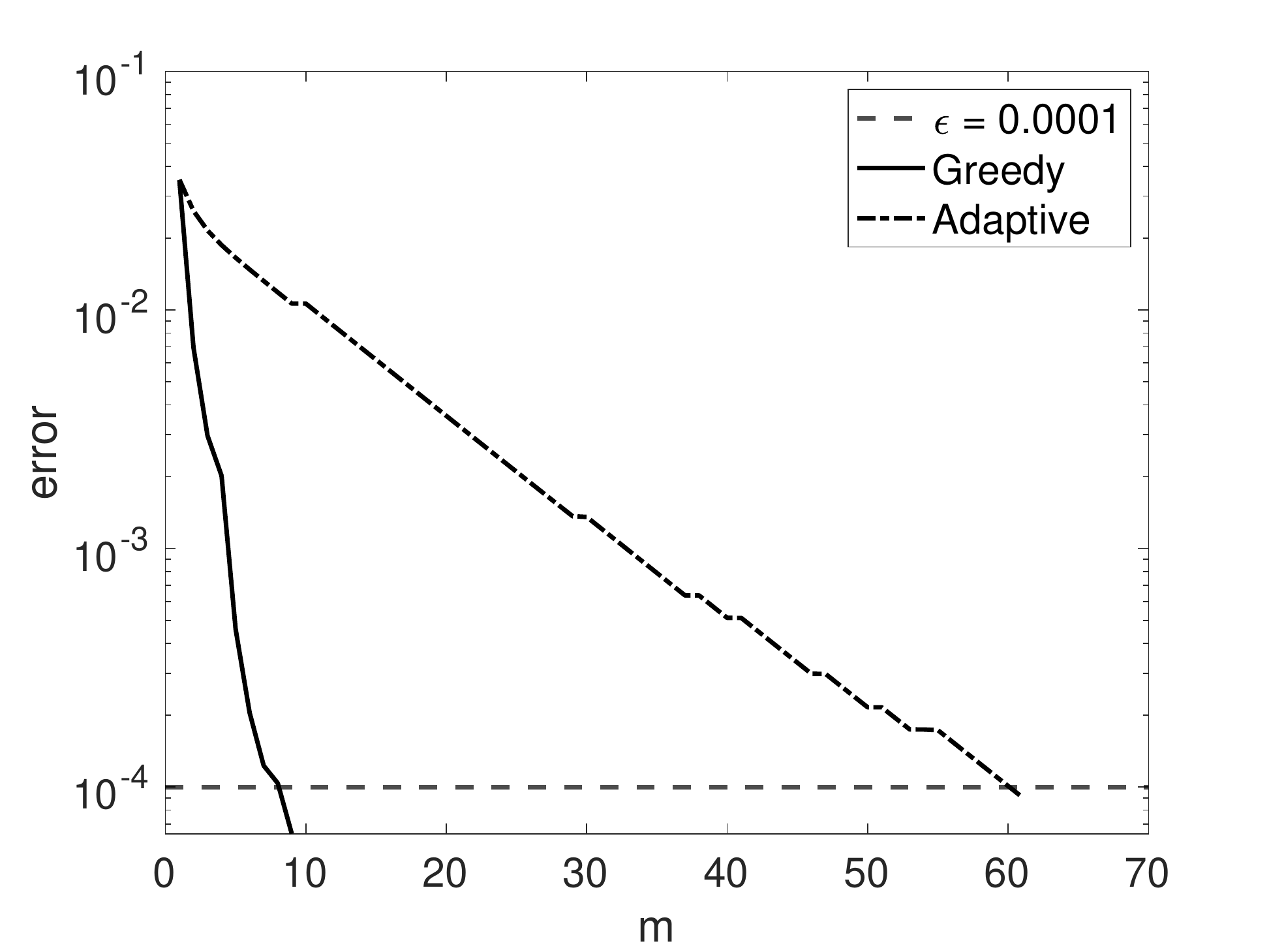}
	\hspace*{-0.5cm}\includegraphics [width=.35\textwidth]{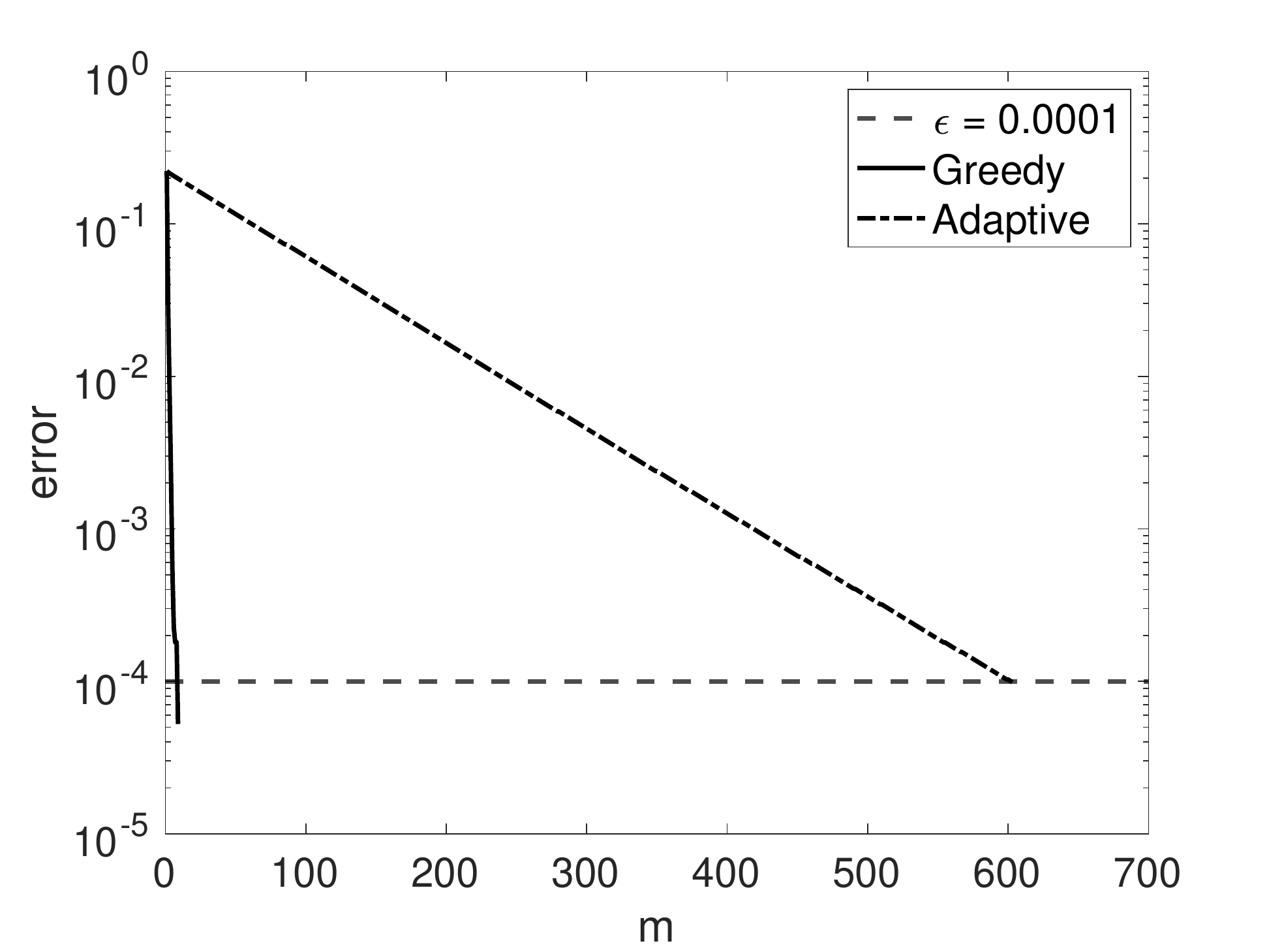}
	\caption{Error between $\bar u$ and $\bar u_n$ versus $n$ for both the (\emph{adaptive}) Taylor and greedy reduced models.  Left: $s=2$, $a_{\rm min}=0.1$; middle: $s=4$, $a_{\rm min}=0.1$; right: $s=4$, $a_{\rm min}=0.01$.} \label{fig:adapt_greedy}
\end{figure}

We next discuss greedy basis constructions. In this case, the reduced linear space $V_n(G)$ is constructed by starting with the function $\vphi_0:=u(0)$ and then use a particular random weak greedy algorithm \footnote{We use a version of the probabilistic weak greedy algorithm given in \cite{CDD}.} to generate the reduced basis { functions $\vphi_1,\dots,\vphi_{n-1}$.  Each $\vphi_j$ is a snapshot $\vphi_j=u(y^{(j)})$ of the solution at a judiciously chosen point $y^{(j)}\in Y$. We denote by $\bar \vphi_j$ an accurate finite element approximation of $\vphi_j$, $j=0,1,\ldots,n-1$, and we define $V_n(G):={\rm span}\{\bar \vphi_0, \bar \vphi_1,\dots,\bar \vphi_{n-1}\}$. The reduced model is now
\be
\label{greedyrm}
 \bar u_n(y):= P_{V_n}(u(y)).
\ee
where $P_{V_n}$ is the Galerkin projection onto $V_n$(G),  namely for a given $y\in Y$, $\bar u_n(y)\in V_n(G)$ is the solution of
\begin{equation*}
\int_D a(\cdot,y)\nabla \bar u_n(y)\cdot\nabla \bar v_n=\int_D f \bar v_n, \quad \bar v_n\in V_n(G).
\end{equation*}
 We compute the error for approximating $u(y)$ using random samples of the parameter $y$ in a similar manner to the Taylor case already discussed.
 
Figure \ref{fig:adapt_greedy} gives a comparison of the performance of} the greedy basis and the (adaptive) Taylor for different values of $s$ and $a_{\rm min}$.
This graph shows that the greedy basis produces a much more accurate reduced model than the Taylor basis given the same allocation $n$ for the dimension of the reduced space.

\subsection{Nonlinear models based on piecewise polynomials}
\label{SS:numericalpp}

The next set of experiments numerically implements a strategy for generating a nonlinear reduced model based on piecewise polynomials similar to that described in \S \ref{sec:pp_PPDE}. We consider the same diffusion coefficients as above and the same values of $s$ and $a_{\rm min}$. 
We again fix a target accuracy $\e=10^{-4}$,
and a target value of $m$ for the dimension of the polynomial space on each cell of the partition. We will see that it is not always possible to achieve a partition of reasonable size if $m$ is chosen to be too small.
This is heuristically clear from the entropy considerations provided in \S \ref{sec:general} and \S\ref{sec:examples}.

Our strategy for generating the partitioning of $Y$ into cells is motivated by the theoretical results of \S \ref{sec:pp_PPDE}. However, we make some modifications of this strategy which we now explain.  Since in our numerical
examples $u$ has singularity near $y=-1$  because  $c_j>0$ for all $j=1,\ldots,64$, we now grade the partition to be finer near $-1$ when we refine a coordinate direction. This is in contrast to the theoretical description, which partitions in a symmetric way for each coordinate $y_j$.  

 On the other hand, we have found that prescribing $\e$ and $m$ and then implementing the theoretical partitioning strategy actually produces a partition with much better accuracy than $\e$,
and thus we have used too many cells. So instead of viewing the target error and $m$ as the parameters to determine the partition, we 
introduce a single parameter $\eta$ to generate a partition.
We then select $\eta$ to give the required accuracy $\e=10^{-4}$ and a good control on $m$ and the number of cells $N$. To be precise, 
 we take $q=1$ and given $\eta>0$, we generate a partition as follows.
\vskip .1in
\noindent
{\bf  Construction of the partition for a given $\eta$ and a non-decreasing sequence $(\rho_j)_{j\geq 1}$:} 
\begin{itemize}
\item[] Choose $J \geq 0$ as the smallest integer such that $\sum_{j=J+1}^{64} \rho_j^{-1}\leq \frac{1}{2} \eta$ and set $\sigma=\frac{\eta}{2J}$, $j=1$;
\item[] While $\sigma \rho_j < 1$ do
\begin{itemize}
\item[] $y_j^0=\frac{1-\sigma\rho_j}{1+\sigma}$, $\lambda_j^0=\sigma(\rho_j +y_j^0)$,  $i =0$;
\item[] While  $y_j^i-\lambda_j^i > -1$
\begin{itemize}
 \item[] Increment $i$;
\item[] Compute $\lambda_j^i=\frac{\sigma}{1+\sigma}(\rho_j+y_j^{i-1}-\lambda_j^{i-1})= \frac{1-\sigma}{1+\sigma}\lambda_j^{i-1}$ and $y_j^i=y_j^{i-1}-\lambda_j^{i-1}-\lambda_j^i$;
\end{itemize}
\item[] End do
 \item[] If $y_j^i-\lambda_j^i < -1$ set $\lambda_j^i=\frac{1}{2}\left(y_j^{i-1}-\lambda_j^{i-1}+1\right)$ and $y_j^i=\frac{1}{2}\left(y_j^{i-1}-\lambda_j^{i-1}-1\right)$;
\item[] Increment $j$;
\end{itemize}
\item[] End do
\item[]  Set $y_l^0 = 0$, $\lambda_l^0 = 1$ for $l=j,...,J$.
\end{itemize}

 The algorithm generates a tensor product partition with cells $Q_\lambda(\bar y)$ of the form~\eqref{def:Q}. For each cell 
$Q_\lambda(\bar y)$ from this partition we define a sequence $(\tilde\rho_j)_{j\geq 1}$, where
\begin{equation*}
\tilde \rho_j:=\left\lbrace
\begin{array}{ll}
\frac{\rho_j+\bar y_j}{\lambda_j}, & \textrm{when } \sigma \rho_j < 1,\\ 
\rho_j, & \textrm{otherwise}.
\end{array}\right.
\end{equation*}
It is easy to check that conditions similar to those in Corollary~\ref{cor1} are satisfied. Namely, $\tilde \rho_j \geq \kappa$, $j=1,\ldots,64$, and $\|(\tilde \rho_j^{-1})_{j =1}^{64}\|_{\ell_q}\leq \|(\rho_j^{-1})_{j =1}^{64}\|_{\ell_q}$. Moreover, we have
$$
 \tilde\delta := \max_{j=1,\ldots,64}\left| \frac{\rho_jc_j+\bar y_jc_j}{1+\bar y_jc_j}\right|<1,
$$
since $\rho_jc_j=1-a_{\rm min}/2<1$, but not necessarily that $\tilde \delta \leq \delta$. However, we can still get the error bound~\eqref{cor11} of Corollary~\ref{cor1}, but with constant $C(\delta, \rho, q)$
replaced by the potentially larger constant $C(\tilde \delta, \rho, q)$. A uniform error bound can be obtained by taking the constant associated to the largest $\tilde \delta$ over all cells in the partition.

Table~\ref{table1} shows the number of terms $m$ needed in the Taylor expansion on each of the $N$ cells from our partition  to meet our error criteria. 
\begin{table}[htbp]
	\centering
	\begin{tabular}{|l|c|c|c||l|c|c|c|}
	\hline
	\multicolumn{4}{|c||}{$a_{min} = 0.1$} & \multicolumn{4}{c|}{$a_{min} = 0.01$} \\
	$\#$ of cells & $s = 2$ & $s = 3$ & $s = 4$ & $\#$ of cells & $s = 2$ & $s = 3$ & $s = 4$ \\
	\hline
	$N = 1$ & 102 & 68 & 61 & $N = 1$ & 666 & 614 & 603 \\
	 $N = 3$ & 29 & 13 &  9 & $N = 3$ & 48& 30 & 27 \\
	$N = 8$ & 22 & 8 & 5 & $N = 10$ & 24 & 11 & 8 \\
	\hline
	\end{tabular}
\caption{ Number of terms $m$ needed to meet the target accuracy $\e=10^{-4}$  on each cell using the piecewise  (\emph{adaptive}) Taylor polynomial approximations.} \label{table1}
\end{table}
 We see that allowing partitioning can significantly reduce the number $m$ of polynomial terms needed to meet the target accuracy. For example,  in the case $N=1$ (i.e., no partitioning), we need to take $m=603$
 whereas using only ten cells the necessary $m$ is reduced to eight.
Note however, that reducing $m$ even further may cause a considerable growth in the number of cells $N$. Finally, we mention that $J=1$ in all the examples above.

\begin{remark}
In the above numerical examples, we have not considered the case of using nonlinear models based on piecewise greedy bases. The reason for this is that we do not have an a priori way to generate a good partition of $Y$ into cells when greedy bases rather than polynomial bases are used on each cell. An appropriate strategy would seem to be to do the partitioning in tandem with the local greedy constructions. Strategies for doing this are currently under investigation.
\end{remark}

\subsection{State estimation using linear and nonlinear reduced models}
\label{SS:numericaldataassim}

As remarked in the introduction, we anticipate that one of the major advantages of using library approximation occurs in the problem of state estimation from data observations. In this section, we recall  the state estimation problem and execute several
numerical experiments indicating the performance of piecewise polynomial approximations for this problem.
	
In state estimation, we are given measurements of an unknown state $u(y^*)$ where $u$ is the solution to 
\eqref{elliptic}
with the model $a$ for the diffusion known to us. We assume that the data is of the form
\be
\nonumber
w_j= l_j(u(y^*)),\quad j=1,\dots,L,
\ee
where the $l_j$ are linear functionals defined on $V$. Each linear functional $l_j$ has a Riesz representation
\be
\nonumber
l_j(v)=\langle v,\omega_j\rangle_V,\quad j=1,\dots,L.
\ee
The functions $\omega_j$, $j=1,\dots,L$, span a subspace $W$ of $V$. Without loss of generality, we can assume that the dimension of $W$ is $L$ since otherwise there is redundancy in the measurements.
	
We want to use these data observations together with the known model $a$ for diffusion in order to construct an approximation $\hat u$ to the state $u(y^*)$. Note that $y^*$ and $u(y^*)$ are not necessarily uniquely determined by the measurements. One way of proceeding, as was proposed in \cite{MPPT}, is to employ a reduced model based on a linear space $V_n$ to approximate $\cM$.
The algorithm in \cite{MPPT} constructs an approximation $\hat u_n$ to $u(y^*)$ by solving a least squares fit to the data from $V_n$. This algorithm was shown to be optimal in a certain sense (see \cite{BCDDPW,DPW}) once $V_n$ is chosen. The performance of this algorithm is upper bounded by
\be
\label{dae}
\|u{ (y^*)}- \hat u_n\|_V\le \mu_n\e_n, \quad {\rm where} \ \e_n:=\dist(\cM,V_n)_V.
\ee
Here $\e_n:=\dist(\cM,V_n)_V$ and $\mu_n=\mu(W,V_n)\ge 1 $ is a certain inf-sup constant which can be interpreted as the reciprocal of the angle between $V_n$ and the space $W$ \cite{BCDDPW17}, namely 
\begin{equation*}
\mu_n=\mu(W,V_n):=\left(\inf_{v\in V_n}\sup_{w\in W}\frac{\langle v,w\rangle_V}{\|v\|_V\|w\|_V}\right)^{-1}.
\end{equation*}	
This motivates choosing a nested sequence $V_1\subset V_2\subset \cdots$ of spaces with $\dim(V_j)=j$ and selecting a space from this sequence
 which minimizes the right side of \eqref{dae}. Note that while $\e_n$ decreases when increasing $n$, the constant $\mu_n$ increases and is in fact infinite if $n>L$.  
	 
For our numerical experiments in state estimation we use the same models for the diffusion $a$ as described in \eqref{affine1}-\eqref{rhoj}.
For the measurements, we take linear functionals which emulate point evaluation. Specifically, each $l_j$ is of the form
\be
\label{ljform}
l_j(u):=\int_{D} u(x)K(x-x_j)\, dx, \quad K(x):= \exp({-\lambda |x|^2}), 
\ee
where $|x|$ is the Euclidean norm of $x$ and $\lambda=227.\bar{5}$.

In our numerical experiments, we set $y^*=0.5384$, but of course operate as if $y^*$ is unknown to us. We take $L=20$ measurements of the form \eqref{ljform}, where the centers $x_j$ are chosen at random, applied to the solution $u(\cdot ,y^*)$ of \eqref{elliptic} with $a$ satisfying \eqref{affine1}-\eqref{rhoj} with $s=4$ and $a_{\rm min}=0.1$. We only see these measurements and not the entire function $u(\cdot,y^*)$.

\begin{figure}[htbp]
	\begin{minipage}{0.5\textwidth}
		\centering
		\hspace*{-1cm}\includegraphics [width=0.8\textwidth]{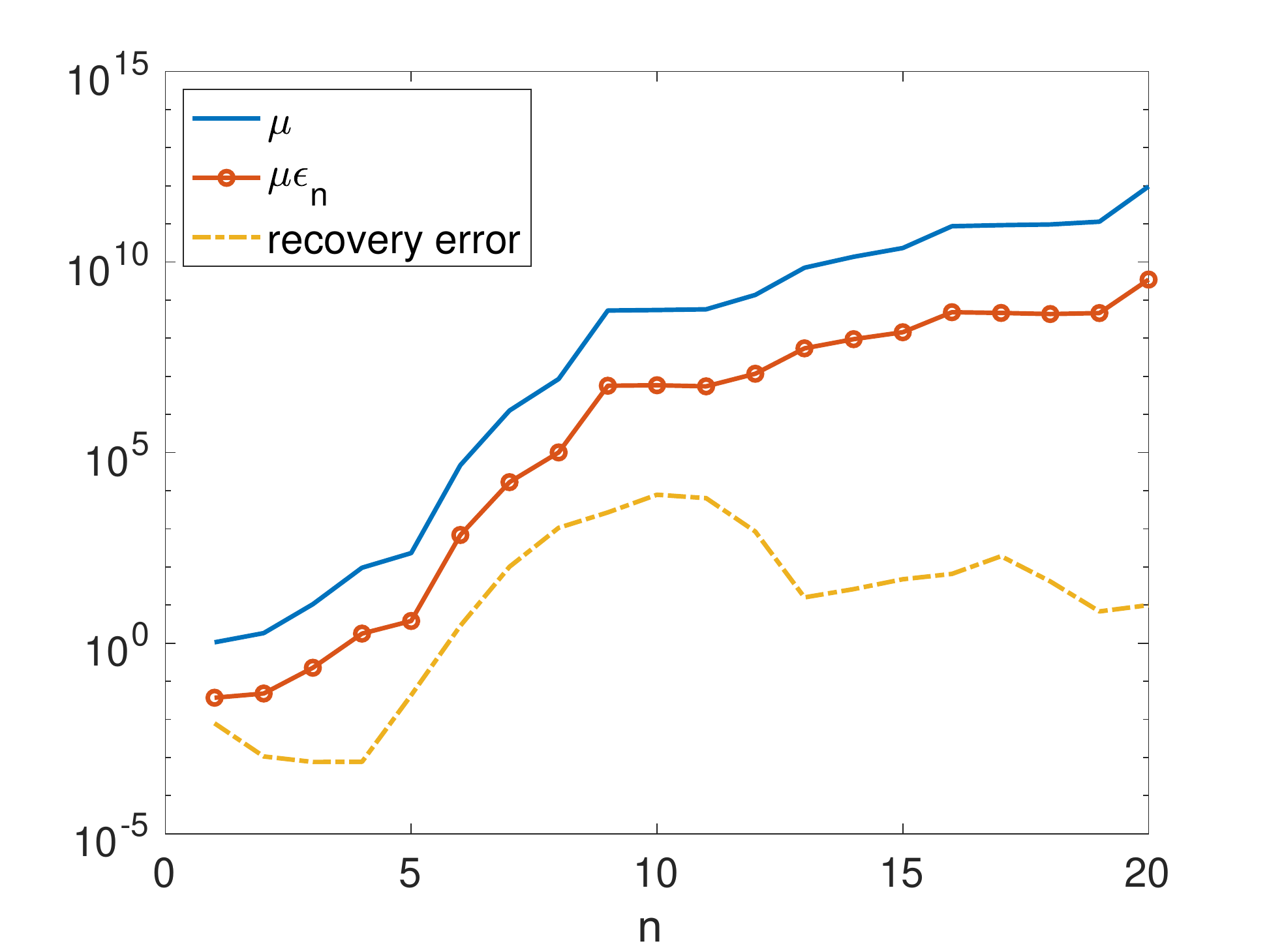}
	\end{minipage}
	\begin{minipage}{.5\textwidth}
		\centering
		\hspace*{-1.5cm}\begin{tabular}{|c|c|c|c|}
		\hline
		$n$ & $\mu_n$ & $\mu_n\epsilon_n$ & $\|\bar u(y^*)-\hat u_n\|_V$ \\
		\hline
		5 & $2.30600\times 10^2$ & $3.81786\times 10^{0}$ & $4.30783\times 10^{-2}$ \\
		10 & $5.48266\times 10^8$ & $5.82909\times 10^{6}$ & $7.82917\times 10^{3}$ \\
		15 & $2.31818\times 10^{10}$ & $1.43338\times 10^{8}$ & $4.74617\times 10^{1}$ \\
		\hline
	\end{tabular}
	\end{minipage}%
	\caption{The constant $\mu_n$, the upper bound $\mu_n\epsilon_n$ and the \emph{recovery error} $\|\bar u(y^*)-\hat u_n\|_V$  for the (\emph{adaptive}) Taylor approximation when $L=20$, $s=4$ and $a_{\rm min}=0.1$. Left: graphs for $n=1,2\ldots,20$; right: values for $n=5,10,15$.} \label{fig:StateEst} 
\end{figure}

Our first numerical experiment is to compute the behavior of $\mu_n$, the recovery error $\|\bar u(y^*)-\hat u_n\|_V$ and its upper bound $\mu_n\epsilon_n$, see \eqref{dae}, for different choices of $V_n$, where $V_n$ is the (adaptive) Taylor with $n$ terms and $\epsilon_n$ is the approximation error computed as discussed in \S\ref{SS:numericallinear}. The values obtained for $n=1,2,\ldots,20$ when $L=20$, $s=4$ and $a_{\rm min}=0.1$ are provided in Figure \ref{fig:StateEst}. The important thing to observe in this figure is that increasing the value of $n$ (in order to improve the approximation error) causes $\mu$ to increase greatly and thereby limiting the recovery accuracy. We shall see in the next experiments that this can be circumvented by using piecewise polynomial approximations.

\begin{figure}[htbp]
	\centering
	\hspace*{-0.1cm}\includegraphics [width=.35\textwidth]{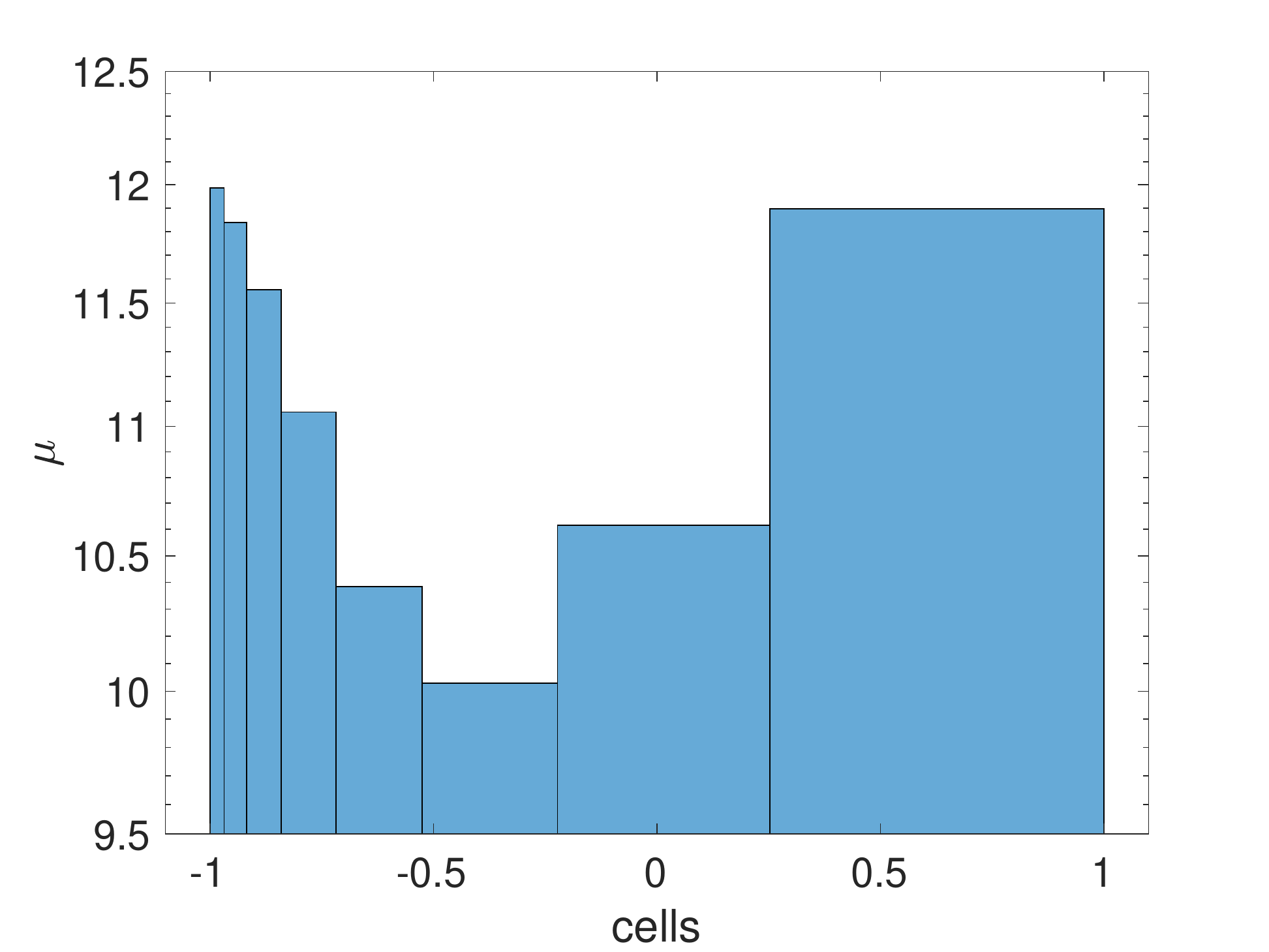}
	\hspace*{-0.5cm}\includegraphics [width=.35\textwidth]{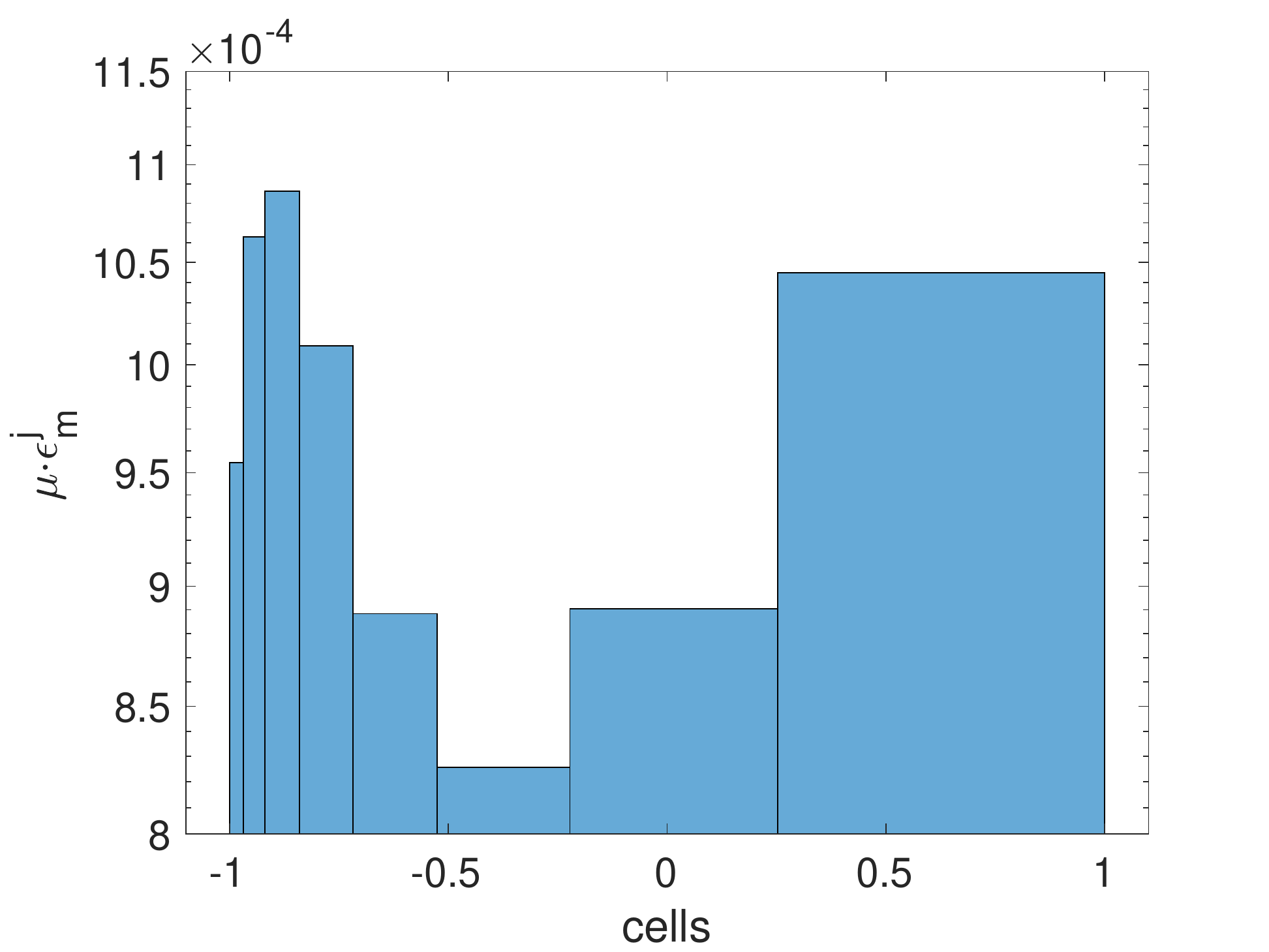}
	\hspace*{-0.5cm}\includegraphics [width=.35\textwidth]{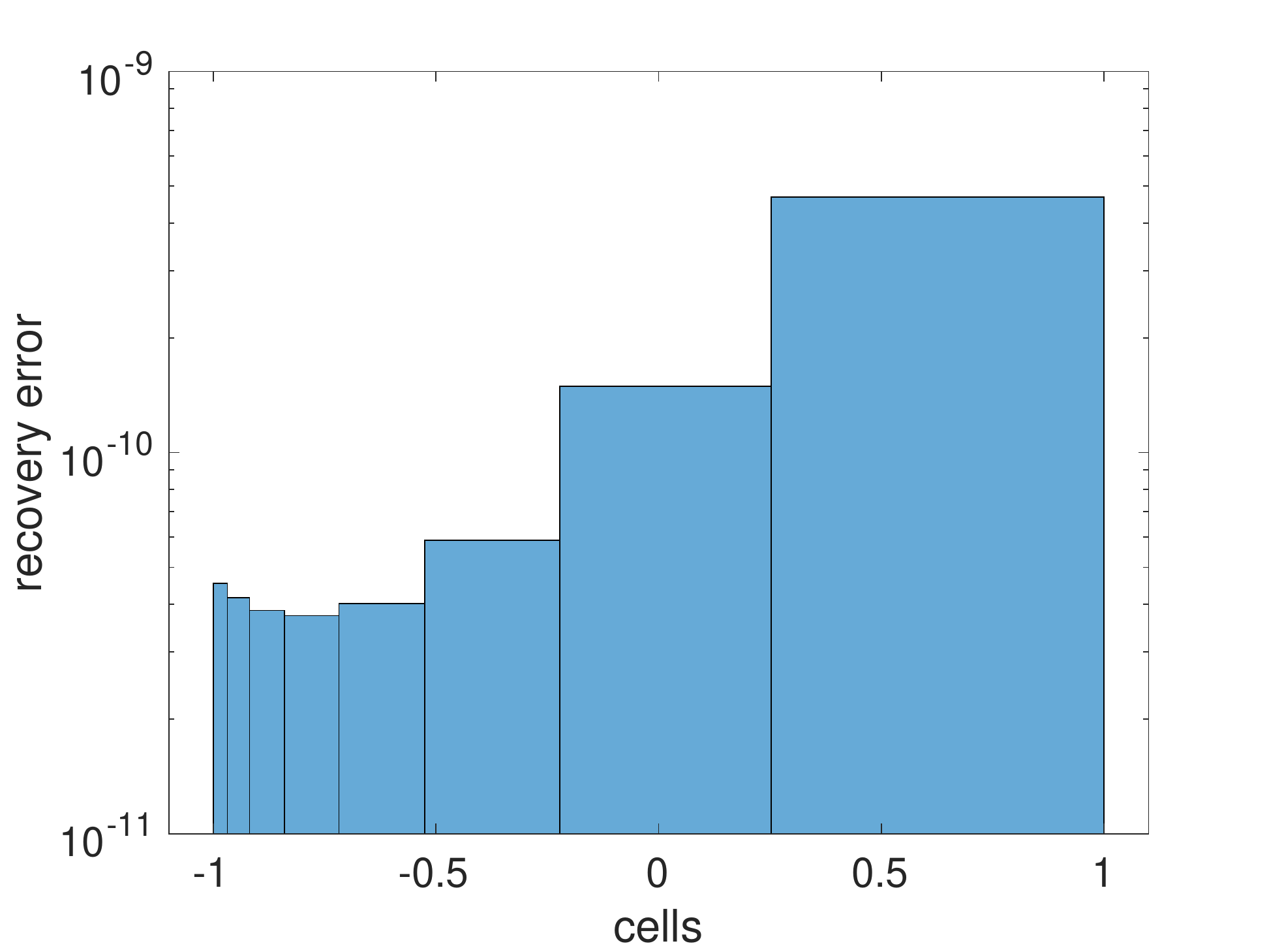}
	\caption{Results of the piecewise (\emph{adaptive}) Taylor polynomial approximation on each cell when $L=20$, $s=4$ and $a_{\rm min}=0.1$. Left: the constant $\mu_m^j$; middle: $\mu_m^j\epsilon_m^j$; right: {\it recovery error} $\|\bar u(y^*)-\hat u_m^j\|_V$.} \label{fig:mu_cells}
\end{figure}

Notice that the dimension $n$ of $V_n$ is limited by $n\le L$ since otherwise $\mu_n$ is infinite. This motivates the use of library approximation with the spaces in the library of small dimension $m\le L$. We do such a numerical experiment using piecewise Taylor polynomial approximation obtained via the adaptive method. We partition $Y$ into 8 cells. This partition corresponds to only subdividing the first coordinate direction $y_1$. Each cell gives rise to a ``local'' value of the inf-sup constant $\mu_m^j:=\mu(W,V_m^j)$, $j=1,2,...,8$, where the $V_m^j$'s are the spaces in the library associated with the partition of $Y$. Finally, we use $m=5$ which ensures that the local approximation error satisfies $\epsilon_m^j\le \e=10^{-4}$ for $j=1,2,...,8$, see Table \ref{table1}. Figure \ref{fig:mu_cells} gives the value of $\mu_m^j$, the upper bound $\mu_m^j\epsilon_m^j$ and the {\it recovery error} $\|\bar u(y^*)-\hat u_m^j\|_V$, $\hat u_m^j\in V_m^j$, for each cell $j=1,2,...,8$. Notice that the values of $\mu$ do not depend on $y^*$.
Also note that the ``local'' constant $\mu$ for the various cells does not exceed $12$ while it was about $230$ for the one cell case, see Figure \ref{fig:StateEst}-right. Moreover, we observe that for all the cells, the upper bound $\mu_m^j\epsilon_m^j$ is smaller than $1.1\times 10^{-3}$, which ensures that the \emph{recovery error} (unknown in practice) is less than $1.1\times 10^{-3}$. Note however that we are not providing an algorithm for determining to which cell the parameter $y^*$ is most likely belong to.
		
	\section{Conclusions}
	\label{S:conclusions} 
	In this section, we briefly discuss the possible advantages and disadvantages of using nonlinear reduced models in the context of parametric PDEs. We consider only the case of elliptic
	PDEs \eqref{elliptic} with affine diffusion coefficients  \eqref{affine}. We suppose that for the given $(\psi_j)_{j\ge1}$, there is a nondecreasing sequence $(\rho_j)_{j\ge 1}$ with $\rho_1>1 $ satisfying \eqref{eqn:WUEA}.	Quantitative theorems for constructing online solvers with performance guarantees are proven using assumptions on the growth of the sequence $(\rho_j)_{j\ge 1}$. A typical assumption that gives a performance guarantee is that the sequence $ (\rho^{-1}_j)_{j\ge 1}$ is in $ \ell_q(\mathbb{N})$ for some $q<2$ (see \cite{BCM,BDGJP2019}). We assume that we have such a sequence with a fixed value of $q$.
	Our discussion is guided by both the theoretical and numerical results of this paper.

\subsection{Offline cost for constructing the solver for linear reduced models} \label{SS:offline}
	
Let us first consider the case where our interest is to construct an online solver for the parametric PDE which performs with a guaranteed approximation error $\e$. There is a distinction in the offline cost of constructing such a solver, depending on whether it is based on Taylor expansions  or on a greedy basis expansion. 
	
When using a Taylor polynomial approximation, we need to find a lower set $\Lambda=\Lambda(\e)$ of indices used in the Taylor polynomial expansion \eqref{ps}, where the approximation to the solution belongs to the space $\cP_\Lambda$. Recall that we presented two methods for finding such an index set $\Lambda$, which we referred to as the \emph{a priori} and the \emph{adaptive} method. The \emph{a priori} method is numerically cheap since it only requires us to sort the $\rho^{-\nu}$ to identify the largest of these numbers (see \cite{BDGJP2019} for one such sorting algorithm). Once the set $\Lambda$ is identified, the Taylor coefficients $\bar t_\nu$ can be computed recursively with finite element solvers as already discussed. The \emph{adaptive} method to build the set $\Lambda$ may seem more expensive as it requires the computation of all $\bar t_\nu$ in the reduced margin of the adaptively constructed monotone set, while only a few may be included in the set 
 $\Lambda_n^*$; compare for instance (7.104) and (7.105) in \cite{CD}. However, this algorithm is preferred in our numerical experiments presented because it generates sets $\Lambda$ with eventually smaller cardinality by assessing precisely the magnitude of $\| \bar t_\nu \|_{V}$ instead of using its upper bound $C_u\rho^{-\nu}$ (see \cite[Lemma 3.14]{CD}). 
		 
Consider next the linear reduced model based on the Galerkin projection onto a linear space $V_{n}$ of dimension $n$ constructed by a weak greedy selection of snapshots from the solution manifold. The advantage of such a greedy construction is that $n$ may be much smaller than the number of terms $\#\Lambda$ used in the Taylor polynomial approximation (see Figure~\ref{fig:adapt_greedy}). Yet, the deficiencies in such greedy algorithms are that the offline cost for the selection of the greedy basis using an $\e$-net training set grows like $O(\e^{-c/r}e^{C\e^{-1/r}})$ (see for instance (8.89) together with (8.108) from \cite{CD}) which may be prohibitive for small $\e$. This of course was one of the main motivations for using nonlinear models in place of linear models.

	\subsection{Offline cost for constructing a solver using nonlinear reduced models }
	\label{SS:offlinenl}

 We discuss next the offline cost in the construction of nonlinear reduced models. Let us first consider reduced models based on piecewise Taylor polynomials. 
 We have given a priori recipes for the tensor product partitioning of $Y$ into cells $Q$ based on the knowledge of the sequence $(\rho_j)_{j\ge 1}$, and thus the main issue is building the appropriate basis 
 for each cell $Q$ of this partition.  
 This requires the computation of the finite element approximation of the appropriate Taylor coefficients on each cell. Note that these computations can be done in parallel. The total cost of this offline construction is governed by the total number $N$ of cells in the partition and the number of terms $m$ used on each cell. In our numerical examples, these constructions were not an issue because the number of cells $N$ was reasonable for 
  moderate values of $m$.

 We have given a priori bounds on the number of cells needed for the partition in \S \ref{SS:apriori}. Recall that if we are in a situation where linear methods (such as polynomial or greedy) give an approximation rate $Mn^{-r}$ then we can guarantee
 an approximation error $\e=n^{-r}$ by using piecewise polynomials with $m$ terms and $N\le e^{C(n/m)^\alpha}$ cells. If we think of the cost of creating a polynomial approximation with $m$ terms to scale like $e^{cm}$, which we know is the case for greedy constructions,  then the cost for constructing the piecewise polynomial is bounded by $e^{C(n/m)^\alpha+cm}$. By choosing $m<n$ appropriately, this is always less than the cost of the approximation without partitioning, which is $e^{Cn}$. For example, if $\alpha=1$ then we could choose $m=\sqrt{n}$ and get the total piecewise polynomial cost to be $e^{C\sqrt{n}}$ as compared with the $e^{Cn}$ if we do not partition. In our numerical
 examples, we have seen that the a priori bounds on the number of cells is quite pessimistic, and we actually get better performance than that predicted by the a priori estimates for the number of cells.

\subsection{Online cost for constructing the approximate solution for linear reduced models} \label{SS:online}

	 If we use a linear reduced model based on Taylor polynomials, then once the index set $\Lambda$ is found and the Taylor coefficients $\bar t_\nu$, $\nu\in\Lambda$, are computed, the reduced model is
	\be
	\nonumber
	\bar u(y)=\sum_{\nu\in\Lambda} \bar t_\nu y^\nu.
	\ee
	Thus, given a parameter query, the online cost for the evaluation of $\bar u(y)$ is trivial.

	If in place of a Taylor polynomial space for the reduced model, we use a greedily generated linear space $V$ of dimension $n$ there are additional online costs.  Given a parameter query $y$ one must find the Galerkin 
	projection of $u(y)$ onto $V$. This entails the inversion of an $n\times n$ dense matrix where the matrix depends on $y$. In certain cases, such as when the diffusion coefficient is affine, this can be somewhat mitigated
	by precomputing certain matrices (see the discussion in \cite{CD}). Therefore, there is a balancing between having a smaller 
	dimensional reduced model (when compared with the polynomial case) and the additional cost of matrix inversion in an online solver. 
	
	Notice also that the accuracy of the online performance given above for reduced models using Taylor polynomials
	can be improved by using a Galerkin projection onto the polynomial space in place of the plug in formula. However, this projection would also involve an expensive matrix inversion.

\subsection{Online cost for constructing the approximate solution for nonlinear reduced models}
 \label{SS:onlinenl}
	
Building an online solver based on piecewise Taylor polynomial approximations proceeds by building a linear solver for each cell of the partition. An additional step is required to determine which space from the library of spaces should be used for the query $y$. This only requires the identification of the cell which contains $y$, and is easily determined from the knowledge of the partition since the cells are hyperrectangles.

 \subsection{Storage costs}
 \label{SS:storage} 

The storage cost for the online solver is dominated by the storage of the basis functions. They are typically large vectors depending on $\e$, $D$ and $f$ in \eqref{elliptic}. We observe from our numerical experiments that the storage cost is higher for linear reduced models using Taylor polynomials compared to the greedy reduced basis algorithm; see Figure~\ref{fig:adapt_greedy}. Moreover, the costs for Taylor polynomial reduced models and piecewise Taylor polynomial reduced models are quite comparable. For example, from Table~\ref{fig:adapt_greedy} we realize that for a target accuracy $\e=10^{-4}$ and $s=3$, $a_{\rm min}=0.01$, the linear reduced model uses 614 basis functions $\bar t_\nu$ while the piecewise Taylor construction has 48 cells with $m=9$ terms on each cell, and hence requires the storage of 432 vectors.

\subsection{Summary}
\label{SS:summary}

 The advantages of  a Taylor polynomial based linear reduced model are:
\begin{itemize}
	\item  possible simple identification of the set $\Lambda$ with no need for optimization or search algorithms;
	\item fast computation of the online solver $\bar u(y)$. 
\end{itemize}
	The deficiency in such constructions is that to reach a small target accuracy $\e$ the dimension $m=\#\Lambda$ may be very large and thus affect the offline construction.
	A large value of $m$ would also affect storage costs.
	
	The advantage of a greedily chosen linear reduced model is that the dimension required for it to reach a target accuracy is typically much smaller than what is required when using Taylor polynomials.
	The disadvantage is the large offline cost to construct the greedy basis when the required dimension is large, along with the higher cost of executing an online solver. There is, however, a 
	savings in storage because the dimension of the greedy space is small.

A piecewise polynomial nonlinear reduced model has the advantage of being able to achieve a better accuracy than linear reduced models while still taking $m$ small, provided that the number of cells $N$ in the piecewise
construction is moderate. In this paper, we have given both a priori bounds on the necessary size of $N$ as well as numerical bounds. Both bounds show the advantage of this approach. The potential deficiency of
this approach is a large storage cost if $N$ is large. Our numerical examples suggest that $N$ is considerably smaller than the a priori bounds thereby making this a viable approach when the desired accuracy $\e$ is
small.

\bibliographystyle{plain}
\bibliography{bibliography}
 
\end{document}